\documentclass[12pt]{amsart}
\usepackage{amssymb}
\newcommand{\zwei}[2]{\left[
    \begin{array}{c} #1 \\ #2
    \end{array} \right]} 

\newcommand{\DOT}{{\setlength{\unitlength}{1pt}\begin{picture}(3.5,2)(1,1)
\put(2.3,2){\circle*{2}}\end{picture}}}

\newcommand{\Fdot}{{F\!_{\DOT}}}

\newcommand{\Kdot}{{K_{\DOT}}}
\newcommand{\Kpdot}{{K_{\DOT}^\perp}}

\newcommand{\TP}{{\mathcal T}{\mathcal P}}

\newcommand{\QED}{
\setlength{\unitlength}{1.0pt}%
\begin{picture}(7.5,7.5)
\put(0,-5){\rule{2.5pt}{2.5pt}}
\put(0,-2.5){\rule{2.5pt}{2.5pt}}
\put(0,0){\rule{5pt}{2.5pt}}
\put(0,2.5){\rule{7.5pt}{2.5pt}}
\end{picture}\vspace{10pt}}

\newtheorem{thm}{Theorem}[section]
\newtheorem{prop}[thm]{Proposition}
\newtheorem{lem}[thm]{Lemma}

\newtheorem{conj}[thm]{Conjecture}
\newtheorem{example}[thm]{Example}
\newtheorem{conex}[thm]{Counterexample}

\newenvironment{ex}{\begin{example}\rm}{\end{example}}
\newenvironment{cex}{\begin{conex}\rm}{\end{conex}}
\newtheorem{definition}[thm]{Definition}
\newenvironment{defn}{\begin{definition}\sl}{\end{definition}}
\newtheorem{remark}[thm]{Remark}
\newenvironment{rem}{\begin{remark}\rm}{\end{remark}}

\begin{document}

\title[Polynomial systems and a conjecture of 
        Shapiro and Shapiro]{Real Schubert Calculus: Polynomial systems 
	and a conjecture of Shapiro and Shapiro} 
\author{Frank Sottile}
\address{\hskip-\parindent
        Frank Sottile\\
        Department of Mathematics\\
        University of Wisconsin\\
        Van Vleck Hall\\
        480 Lincoln Drive\\
        Madison, Wisconsin 53706-1388\\
        USA}
\email{sottile@math.wisc.edu}
\urladdr{http://www.math.wisc.edu/\~{}sottile}
\thanks{Supported in part by NSERC grant  OGP0170279 \ 
and NSF grants DMS-9022140 and DMS-9701755}

\begin{abstract}
Boris Shapiro and Michael Shapiro have a conjecture 
concerning the Schubert calculus and real enumerative geometry and 
which would give infinitely many families of 
zero-dimensional systems of real polynomials (including families of
overdetermined systems)---all of whose solutions are real.
It has connections to the pole placement problem in linear 
systems theory and to totally positive matrices.
We give compelling computational evidence for its validity, prove it
for infinitely many families of enumerative problems, show how a simple
version implies more 
general versions, and present a counterexample to a 
general version of their conjecture.

This is a companion paper to~\cite{FRZ} and~\cite{Verschelde}, which
describe the mathematics involved in two spectacular computations
verifying specific instances of this conjecture.
\end{abstract}

\maketitle

\section{Introduction}
Determining the number of real solutions to a system of polynomial equations
is a challenging  problem in symbolic and numeric
computation~\cite{G-VRRT,Sturmfels_monthly} with real world applications.
Related questions include when a problem of 
enumerative geometry can have all solutions real~\cite{Sottile97c} and
when may a given physical system be controlled by real output
feedback~\cite{Byrnes,RSW,SADG}.
In May 1995, Boris Shapiro and Michael Shapiro
communicated to the author a remarkable conjecture connecting these three
lines of inquiry. 

They conjectured a relation between topological invariants of the real and
of the complex points in an intersection of Schubert cells in a flag
manifold, if the cells are chosen according to a recipe they give.
When the intersection is
zero-dimensional, this asserts that all points are real.
Their conjecture is false---we give full description and present a
counterexample in Section 5. 
However, there is considerable evidence for their conjecture if the Schubert
cells are in a Grassmann manifold.
It is this variant which is related to the lines of inquiry 
above and which this paper is about.

Here is the simplest (but still very interesting and open) special case of
this conjecture:
Let $m,p>1$ be integers and let $X$ be a $p\times m$-matrix of
indeterminates.
Let $K(s)$ be the $m\times(m+p)$-matrix of polynomials in $s$
whose $i,j$th entry is
\begin{equation}\label{eq:K-matrix}
     \binom{j-i}{i-1} s^{j-i}.
\end{equation}
Set 
$$
\varphi_{m,p}(s;X)\ :=\  \det\zwei{K(s)}{I_p\quad X},
$$
where $I_p$ is the $p\times p$ identity matrix.

\begin{conj}[Shapiro-Shapiro]\label{conj:shapiroI}
For all integers $m,p>1$, the polynomial system 
\begin{equation}\label{eq:conjectureI}
\varphi_{m,p}(1;X)\ =\ \varphi_{m,p}(2;X)\ =\ \cdots\ =\ 
\varphi_{m,p}(mp;X)\ =\ 0
\end{equation}
is zero-dimensional with 
\begin{equation}\label{grassdeg} d_{m,p} \quad := \quad
   \frac{1! \, 2! \, 3! \cdots (p\!- \!2) ! \, (p \!-\!1)! \cdot
       (mp)!}{m!\, (m \! + \! 1)! \, (m \! + \, 2)! 
        \cdots(m \! + \! p \! - \! 1)!} 
\end{equation}
solutions, and all of them are real.
\end{conj}

It is a Theorem of Schubert~\cite{Schubert_degree} that $d_{m,p}$ is a sharp
bound for the number of isolated solutions.
Conjecture~\ref{conj:shapiroI} has been verified for all $1<m\leq p$ with
$mp\leq 12$. 
The case of $(m,p)=(3,4)$ $(d_{m,p}=462)$ is due to an heroic calculation
of Faug\`ere, Rouillier, and Zimmermann~\cite{FRZ} (see
Section~\ref{sec:computation} for a discussion). 
\medskip

Conjecture~\ref{conj:shapiroI} is related to a question of 
Fulton~\cite[\S7.2]{Fulton_introduction_intersection_new}:
``How many solutions to a problem in enumerative geometry
may be real, where that problem consists of counting figures of
some kind having a given position with respect to some given (fixed)
figures.''.
For 2-planes having a given position with respect to fixed
linear subspaces, the answer is that all may be
real~\cite{Sottile97a}.
This was also shown for the problem of 3264 plane conics tangent to five
given conics~\cite{RTV}.
More examples, including that of 3-planes in 
${\mathbb C}^6$ meeting 9 given 3-planes nontrivially, 
are found in~\cite{Sottile97c,Sottile97b}.
The result~\cite{FRZ} extends this to 3-planes in  
${\mathbb C}^7$ meeting 12 given 4-planes nontrivially.

Only the simplest form of the conjecture of Shapiro and Shapiro
has appeared in print~\cite{HSS,RS98,Sottile97c}.
While more general forms have circulated informally, there is no
definitive source describing the 
conjectures or the compelling evidence that has accumulated
(or a counterexample to the original conjecture).
The primary aim of this paper is to rectify this situation and make these
conjectures available to a wider audience.
\medskip

In Section 2, we describe a version of the conjecture 
related to the pole placement problem of linear systems theory.
For this, the integers $1,2,\ldots,mp$ in the polynomial
system~(\ref{eq:conjectureI}) of Conjecture~\ref{conj:shapiroI} are
replaced by generic real numbers and all $d_{m,p}$ solutions are asserted to
be real.
We present evidence (computational and Theorems) in support of it.
Subsequent sections describe the conjecture in greater
generality---for enumerative problems arising from the Schubert calculus on 
Grassmannians in Section~\ref{sec:LR} and a newer extension involving
totally 
positive matrices~\cite{Ando} in Section~\ref{sec:tot_pos}. 
We describe and give evidence for each extension
and show how the version of the conjecture in Section 2 implies more
general versions involving Pieri-type enumerative problems.
In Section~\ref{sec:counter}, we present a counterexample to
their original conjecture and discuss further questions.

A remark on the form of these conjectures is warranted.
Conjecture~\ref{conj:shapiroI} gives an infinite list of specific polynomial
systems, and conjectures that each has only real solutions.
The full conjectures are richer.
For each collection of {\it Schubert data}, Shapiro and Shapiro give a
continuous family of polynomial systems and conjecture that each of the
resulting systems of polynomials has only real solutions.
Conjecture~\ref{conj:shapiroI} concerns one specific polynomial system in
each family, for an infinite subset of Schubert data.

Results here were aided or are due to computations.
Further documentation including Maple V.5 and Singular
1.2.1~\cite{SINGULAR} scripts used are available on the web
page~\cite{Sottile_shap-www}. 

\tableofcontents

\section{Linear equations in Pl\"ucker coordinates}\label{sec:hypersurface}
\subsection{Some enumerative geometry}
Consider the following problem in enumerative geometry:
How many $p$-planes meet $mp$ general $m$-planes in ${\mathbb C}^{m+p}$
nontrivially? 

The set of $p$-planes in ${\mathbb C}^{m+p}$, $\mbox{\it Grass}(p,m+p)$, is
called the {\it Grassmannian of $p$-planes in  ${\mathbb C}^{m+p}$}.
This complex manifold of dimension $mp$ is an algebraic subvariety
of the projective space ${\mathbb P}^{\binom{m+p}{p}-1} $.
To see this, represent a $p$-plane $X$ in ${\mathbb C}^{m+p}$ as the
row space of a $p\times(m+p)$-matrix, also written $X$.
The maximal minors of $X$ are its {\it Pl\"ucker coordinates} and 
determine a point in ${\mathbb P}^{\binom{m+p}{p}-1}$.
This gives the Pl\"ucker embedding of $\mbox{\it Grass}(p,m+p)$.
If $X$ is generic, then its first $p$ columns are linearly
independent, so we may assume they form a $p\times p$-identity matrix.
The remaining $mp$ entries determine $X$ uniquely and give local coordinates
for $\mbox{\it Grass}(p,m+p)$, showing it has dimension $mp$. 

Consider a $m$-plane $K$ to be the row space of a $m\times(m+p)$-matrix,
also written $K$.
Then $K\cap X$ is nontrivial if and only if 
$$
\det\zwei{K}{X}\ =\ 0.
$$
Laplace expansion along $X$
gives a linear equation in the Pl\"ucker coordinates of $X$.

If $K_1,\ldots,K_{mp}$ are $m$-planes in general position, then the
conditions that $X$ meet each of the $K_i$ nontrivially are $mp$
linear equations in the Pl\"ucker coordinates of $X$,
and these are independent by Kleiman's Transversality
Theorem~\cite{Kleiman}. 
Hence there are finitely many 
$p$-planes $X$ which meet each $K_i$ nontrivially
and this number is the
degree of $\mbox{\it Grass}(p,m+p)$ in ${\mathbb P}^{\binom{m+p}{p}-1}$,
which Schubert~\cite{Schubert_degree} determined to be
$d_{m,p}$.

\subsection{The conjecture of Shapiro and Shapiro}\label{sec:shapiroII}
Shapiro and Shapiro gave a recipe for selecting real
$m$-planes $K_1,\ldots,K_{mp}$ and conjecture that when they are in 
in general position, all $d_{m,p}$ $p$-planes meeting each $K_i$ are real.
The standard rational normal curve is the image of the map 
$\gamma:{\mathbb R}\rightarrow {\mathbb R}^{m+p}$
given by 
 \begin{equation}\label{eq:afine_rat}
   \gamma\ :\   s\ \longmapsto\ (1,s,s^2,\ldots,s^{m+p-1})
 \end{equation}
Then the matrix $K(s)$ of the Introduction~(\ref{eq:K-matrix}) has
rows 
$$
\gamma(s),\ \gamma'(s),\ \frac{\gamma''(s)}{2},\ \ldots,\ 
\frac{\gamma^{(m-1)}(s)}{(m-1)!},
$$
where we take derivatives with respect to the parameter $s$.
Let $X$ be a $p\times m$-matrix of
indeterminates.
Define 
$$
\varphi_{m,p}(s;X)\ :=\  \det\zwei{K(s)}{I_p\ X}.
$$

\begin{conj}[Shapiro-Shapiro]\label{conj:shapiroII}
For all integers $m,p>1$ and almost all 
distinct real numbers $s_1,\ldots,s_{mp}$, 
the system of $mp$ equations
\begin{equation}\label{eq:shapiroII}
  \varphi_{m,p}(s_1;X)\ =\ \varphi_{m,p}(s_2;X)\ =\ \cdots\ =\ 
  \varphi_{m,p}(s_{mp};X)\ =\ 0
\end{equation}
is zero-dimensional with $d_{m,p}$ real solutions.
\end{conj}

Let $K(s)$ denote both the $m\times (m+p)$-matrix defined above and its row
space, an $m$-plane. 
Conjecture~\ref{conj:shapiroII} asserts that the $m$-planes
$K(s_1),\ldots,K(s_{mp})$ are in general position, and any $p$-plane
meeting each $K(s_i)$ is real.
The systems are zero-dimensional~\cite{Brockett_Byrnes,EH83}
and there are generically no multiplicities. 
Conjecture~\ref{conj:shapiroI} is the special case when $s_i=i$.

\begin{ex}\label{ex:22}
We establish Conjecture~\ref{conj:shapiroII} when $m=p=2$.
Then 
$$
\varphi_{2,2}(s;X)\ =\ 
\det\left[\begin{array}{cccc}
1&s&s^2&s^3\\
0&1&2s&3s^2\\
1&0&x_{11}&x_{12}\\
0&1&x_{21}&x_{22}\end{array}\right]
$$
is
$$
s^4-2s^3x_{21}+s^2x_{22}-3s^2x_{11}+2s x_{12}+x_{11}x_{22}-x_{12}x_{21}.
$$

We show that if $s,t,u,v\in{\mathbb R}$ are distinct, then the system of
polynomial equations 
\begin{equation}\label{eq:system22}
\varphi_{2,2}(s)\ =\ \varphi_{2,2}(t)\ =\ \varphi_{2,2}(u)\ =\ 
\varphi_{2,2}(v)\ =\ 0
\end{equation}
has all $d_{2,2}=2$  solutions real.
Our method will be to solve~(\ref{eq:system22}) by elimination.

Let $e_i$ be the $i$th elementary symmetric polynomial in $s,t,u,v$.
In the lexicographic term order with $x_{11}>x_{12}>x_{22}>x_{21}$
on the ring ${\mathbb Q}(s,t,u,v)[x_{11},x_{12},x_{22},x_{21}]$, the ideal 
$\langle\varphi_{2,2}(s),\varphi_{2,2}(t),
\varphi_{2,2}(u),\varphi_{2,2}(v)\rangle$
has a Gr\"obner basis consisting of the following polynomials:
$$
2x_{21}-e_1,\  \ 
x_{22}-3x_{11}-e_2,\  \
2x_{12}+e_3,\ \ \mbox{ and }\ \
12x_{11}^2 + 4e_2 x_{11}+ e_1e_3 - 4e_4.
$$
Thus, for distinct $s,t,u,v$,  the system~(\ref{eq:system22})
has 2 solutions and they are real if the discriminant
of the last  equation,
$$
16e_2^2 - 48e_1e_3 + 192 e_4,
$$
is positive.
Expanding this discriminant in the parameters $s,t,u,v$, we obtain
$$
8\left(
(s-t)^2(u-v)^2 + (s-u)^2(t-v)^2 + (s-v)^2(t-u)^2
\right).
$$
Hence all solutions are real, establishing Conjecture~\ref{conj:shapiroII}
when $m=p=2$. 
Theorem~\ref{thm:23} proves Conjecture~\ref{conj:shapiroII}
when $(m,p)=(2,3)$.
\end{ex}

\subsection{Pole placement problem}\label{sec:control}
Suppose we have a physical system (for example, a mechanical linkage) with
inputs 
$u\in {\mathbb R}^m$ and outputs $y\in {\mathbb R}^p$ for which there are
internal states $x\in {\mathbb R}^n$ such that the system evolves 
by the first order linear differential equation
\begin{equation}\label{linsystem}
\begin{array}{rcl}
\dot{x}&=&Ax + Bu,\\ 
    y&=&Cx.
\end{array}
\end{equation}
(We assume $n$ is the minimal number of
internal states needed to obtain a first order equation.)
If the input is controlled by constant output feedback, $u=Xy$,
then we obtain 
$$
\dot{x}\ =\ (A+BXC)x.
$$
The natural frequencies of this controlled system are the roots 
$s_1,\ldots,s_n$ of 
\begin{equation}\label{charpoly}
\varphi(s)\ :=\  \det(sI_n-A-BXC).
\end{equation}
The pole assignment problem asks the inverse question:
Given a system~(\ref{linsystem}) and a
polynomial $\varphi(s)$ of degree $n$, which feedback laws $X$ 
satisfy~(\ref{charpoly})?

A coprime factorization of the transfer function is two matrices
$N(s)$, $D(s)$ of polynomials with $\det(D(s))= \det(sI_n-A)$
and  $N(s)D(s)^{-1}=C(sI_n-A)^{-1}B$.
This always exists.
A standard transformation~({\em cf.}~\cite[\S 2]{Byrnes})
shows that, up to a sign of $\pm 1$, 
 \begin{equation}\label{schubert_form}
   \varphi(s)\ =\ \det\left[
   \begin{array}{cc}  N(s) & D(s)\\ I_p & X \end{array}\right],
 \end{equation}
If we set $K(s)\ :=\ [N(s) \ D(s)]$, write $K(s)$ for the $m$-dimensional
row space of this matrix,  and let $X$ be the 
$p$-plane $[I_p \ X ]$,
then~(\ref{schubert_form}) is equivalent to 
 \begin{equation}\label{eq:control}
    X\cap K(s_i)\ \neq \ \{0\}\quad\mbox{for}\quad i\ =\ 1,\ldots,n,
 \end{equation}
where $s_1,\ldots,s_n$ are the roots of $\varphi(s)$.

If the $m$-planes $K(s_1),\ldots,K(s_n)$ are in general position, then
$mp\geq n$ is necessary for there to be any feedback laws.
These $m$-planes are not {\it a priori} in general position.

To see this, 
let $K:{\mathbb P}^1\rightarrow\mbox{\it Grass}(m,m+p)$ be the
extension of the map  given by
$s\mapsto K(s)$.
Then $K$ is a parameterized rational curve of degree $n$ in
$\mbox{\it Grass}(m,m+p)$.
The space of all such curves $K$ with $n$ distinguished points
$\{K(s_1),\ldots,K(s_n)\}$ has dimension~\cite{Stromme}
$$
mp + n(m+p) + n.
$$
The space of all $n$-tuples of  $m$-planes has dimension $nmp$.
Therefore when 
$$
n \ >\ mp/(mp-m-p-1), 
$$
such $n$-tuples constitute a proper subvariety of 
all $n$-tuples of $m$-planes.

However, the General Position Lemma~\cite{Byrnes} (see
also~\cite{EH83}) states that there is a 
Zariski open subset of the data $A,B,C,\varphi$ such that the $m$-planes
$K(s_1),\ldots,K(s_n)$ are in general position in that the set of  
$X$ satisfying (\ref{eq:control}) has dimension $mp-n$.

Since all rational curves 
$K: {\mathbb P}^1 \rightarrow \mbox{\it Grass}(p,m+p)$ of degree $n$ 
with $K(\infty)=[0\ I_p]$
arise in this way~\cite{Hermann_Martin}, the polynomial
systems of 
Conjecture~\ref{conj:shapiroII} are instances of the pole placement
problem.
Interestingly, these very systems figure prominently in a proof of the General
Position Lemma~\cite{Byrnes_80}.
\smallskip

An important question is whether a given real system may be controlled by
{\em real} feedback~\cite{Byrnes_real,RSW,RS98,SADG,Willems_Hesselink}:
If all roots of $\varphi(s)$ are real, are there any real
feedback laws $X$ satisfying~(\ref{schubert_form})?  
Few specific examples
have been
computed~\cite{BS_homotopy,MWA,RS98,Willems_Hesselink}.  
In~\cite{RS98} an attempt was made to gauge how likely it is for a real
system to be controllable by real feedback and how many of the feedback laws
are real---in the case of $(m,p)=(2,4)$ so that $d_{m,p}=14$.
In all, 600 different curves $K(s)$ were generated, and each of these were
combined with 25 polynomials $\varphi(s)$ having 8 real roots.
Only 7 of the resulting 15,000 systems had all feedback laws real.
This is in striking contrast to the systems given in
Conjecture~\ref{conj:shapiroII}, where all the feedback laws are conjectured  
to be real.

\subsection{Computational evidence}\label{sec:computation}
Consider~(\ref{schubert_form}) as
a map $\mbox{\it Grass}(p,m+p)\rightarrow{\mathbb  P}^{mp}$ in local
coordinates which associates a $p$-plane $X$ to a polynomial $\varphi$ 
(modulo scalars) of degree at most $mp$.
When $K(s)$ is the curve $K_{m,p}(s)$ of Conjecture~\ref{conj:shapiroII},
the inverse image of the polynomial 1 is the single real point $[0\ I_p]$.
Rosenthal suggested that the fibre
over a nearby polynomial may consist of $d_{m,p}$ real points.

Inspired by this, Rosenthal and Sottile~\cite{RS98}
tested and verified several thousand instances of
Conjecture~\ref{conj:shapiroII} when $(m,p)=(2,4)$. 
Each was a specific choice of $m,p$ and $mp$ distinct real numbers 
$s_1,\ldots,s_{mp}$ for which we showed all solutions
to~(\ref{eq:shapiroII}) are real. 
Any verified instance implies 
that all nearby instances in the space of parameters $s_1,\ldots,s_{mp}$
has all of its solutions real. 
In light of the computations described in Section~\ref{sec:control}, we felt  
this provided overwhelming evidence for the validity of
Conjecture~\ref{conj:shapiroII}.

Our method was to solve the
polynomial systems by elimination (see~\cite[\S 2]{CLO_using} for a
discussion of methods to solve systems of polynomial
equations). 
We first choose distinct integral values of the parameters $s_i$ and generate
the resulting system of integral polynomial equations.
Since we are performing an exact symbolic computation, we necessarily work
with integral polynomials.
Next, we compute an eliminant, a univariate polynomial $g(x)$ with the
property that its roots are the set of $x$-coordinates of solutions to our
system.  
When $g(x)$ has $d=d_{m,p}$ roots (Schubert's bound), there is a
lexicographic Gr\"obner basis satisfying 
the Shape Lemma, since this system is zero-dimensional~\cite{EH83}.
It follows that the solutions are rational functions (quotients of integral
polynomials) of the roots of $g(x)$.
In some instances, the eliminant we calculated did not have $d$ roots.
For these we found a different eliminant with $d$ roots.
Lastly, we checked that these eliminants had only real roots.

Table~\ref{table:instances} gives the number of instances we know have 
been checked.
By Lemma~\ref{lem:ordering}(ii), there is a bijection
between instances of $(m,p)=(a,b)$ and $(m,p)=(b,a)$.
Table~\ref{table:instances} also lists the running time to
compute a degree reverse lexicographic Gr\"obner basis for the systems of
Conjecture~\ref{conj:shapiroI}, and the size of that basis.  
This used  Singular-1.2.1~\cite{SINGULAR} on a K6-2-300 processor with 256M
running Linux.
The instances reported in the last 3 columns are not due
the the author.
A more complete account is found in~\cite{Sottile_shap-www}. 
\begin{table}[htb]
 \begin{tabular}{|c||c|c|c|c|c|c|c|}
  \hline$m,p$&4,2&5,2&3,3&6,2&7,2&4,3&2,8\\\hline\hline
  $d_{m,p}$&14&42&42&132&429&462&1430\\\hline
  \# checked&$>12,000$&1000&550&55&2&2&1\\\hline
  time (sec)&.04&1.42&1.50&78.6&8175&--&--\\\hline
  size&1.4K&12.8K&18.6K&202K&4.58M&32M&--\\\hline
 \end{tabular}\bigskip
 \caption{Instances checked\label{table:instances}}
\end{table}

The computations of the last 2 columns stand out.
The first is the case $8,2$ (also one instance each of $7,2$ and $4,3$)
computed by Jan Verschelde~\cite{Verschelde} using his
implementation of the SAGBI homotopy algorithm in~\cite{HSS}.
Since the polynomial system of Conjecture~\ref{conj:shapiroII} was
ill-conditioned, he instead used the equivalent system of
Conjecture~\ref{conj:shapiroII}$'$ (in Section~\ref{sec:simpler} below),
where the $P_i(s)$ were the Chebyshev polynomials. 
These numerical calculations give approximate
solutions whose condition numbers determine a neighborhood
containing a solution.
The solutions of this real system are stable under complex 
conjugation, so it sufficed to check that each neighbourhood and its complex
conjugate were disjoint from all other neighborhoods.
This computation took approximately 25 hours on a 166MHz Pentium II
processor with 64M running Linux.
These algorithms are `embarrassingly parallelizable', and
in principle they can be used to check far larger polynomial systems.

The second is the case of $(m,p)=(3,4)$ of
Conjecture~\ref{conj:shapiroI}
(also all smaller cases with $m\leq p$), computed by Faug\`ere, Rouillier, and
Zimmermann~\cite{FRZ}.   
They first used FGB~\cite{Faugere_FGB} to calculate a degree reverse
lexicographic Gr\"obner basis for the system~(\ref{eq:conjectureI}) for 
$(m,p)=(3,4)$ with $s_i=i$. 
This yielded a Gr\"obner basis of size 32M.
They then computed a rational univariate 
representation~\cite{Rouillier_RUR} (a sophisticated substitute for an
eliminant) in two ways.
Once using a multi-modular implementation of the
FGLM~\cite{FGLM} algorithm and a second time using RS, 
an improvement of the RealSolving
software~\cite{Rouillier_RS} under development.
The eliminant had degree 462 and size 3M, thus its general
coefficient had 2,000 digits.
Using an early implementation of Uspensky's algorithm, they
verified that all 
of its zeroes were real, proving 
Conjecture~\ref{conj:shapiroI} for $(m,p)=(3,4)$. 
In the course of this calculation, they found it necessary
to rewrite their software.

\subsection{Equivalent systems}\label{sec:simpler}
The extension of the map~(\ref{eq:afine_rat}) to ${\mathbb P}^1$
$$
  \gamma\ :\ 
  [t,s]\ \longmapsto\ 
  [t^{m+p-1},\,st^{m+p-2},\,\ldots,\,s^{m+p-2}t,\,s^{m+p-1}]
$$
is a parameterization of the standard real rational normal curve in 
${\mathbb P}^{m+p-1}$
and $K(s)$ is the $m$-plane osculating this 
curve at the point $\gamma[1,s]$.
In general, a parameterized real rational normal curve is a map
$\gamma:{\mathbb P}^1\rightarrow{\mathbb P}^{m+p-1}$ of the form
$$
  [t,s]\ \longmapsto\ [P_1(s,t),\,P_2(s,t),\,\ldots,\,P_{m+p}(s,t)]
$$
where $P_1(s,t),\ldots,P_{m+p}(s,t)$ form a basis for the space of
real homogeneous polynomials in $s,t$ of degree $m+p-1$.
All parameterized real rational normal curves are conjugate by a real
projective transformation of ${\mathbb P}^{m+p-1}$.
We deduce
\medskip

\noindent{\bf Conjecture~\ref{conj:shapiroII} (Geometric form)}
{\it
For all integers  $m,p>1$ and almost all choices of 
$mp$ $m$-planes  $K_1,\ldots,K_{mp}$  osculating
a real rational normal curve at distinct real points, there are exactly
$d_{m,p}$ $p$-planes $X$ satisfying
$$
X\cap K_i\ \neq\  \{0\}\quad\mbox{ for }\quad i=1,\ldots,mp
$$ and all of these $p$-planes $X$ are real.
}\medskip

Thus Conjecture~\ref{conj:shapiroII} is  equivalent to
a conjecture concerning a much richer class of polynomial systems.
\medskip

\noindent{\bf Conjecture~\ref{conj:shapiroII}$'$. }
{\it 
Suppose $m,p>1$ are integers
and $P_1(s),\ldots,P_{m+p}(s)$ are a basis of the space of real polynomials
of degree at most $m+p-1$.
Let $K(s)$ be the $m\times (m+p)$ matrix of polynomials whose $i,j$th entry
is $P^{(i-1)}_j(s)$.
Set
$$
\varphi(s;X)\ :=\ \det\zwei{K(s)}{I_p \ X}.
$$
Then, for almost all choices of distinct real numbers $s_1,\ldots,s_{mp}$,
the system
$$
\varphi(s_1;X)\ =\ \varphi(s_2;X)\ =\ \cdots\ =\ \varphi(s_{mp};X)\ =\ 0
$$
has exactly $d_{m,p}$ solutions, and all of them are real.
}\medskip

The polynomial matrix $K(s)$ of Conjecture~\ref{conj:shapiroII}$'$ differs
from that of Conjecture~\ref{conj:shapiroII} by right multiplication by an
invertible $(m+p)\times(m+p)$-matrix.
Thus the resulting polynomial systems differ primarily by choice of local
coordinates for the Grassmannian.
In linear systems theory, two physical systems are output
feedback-equivalent if their matrices of coprime factors $[N(s)\ D(s)]$
differ in this manner~\cite{RRH}.

We give an equivalent conjecture concerning a simpler system of polynomials 
with 2 fewer equations and unknowns.
We may reparameterize the curve
$K(s)$ of Conjecture~\ref{conj:shapiroII}
and assume $s_{mp-1}=0$ and $s_{mp}=\infty$.
Observe that $K(0)=[I_p\ 0]$ and $K(\infty)=[0\ I_p]$.
The collection of all $p$-planes $X$ satisfying
 \begin{equation}\label{eq:small}
   X\cap[I_p\ 0]\neq\{0\} \quad\mbox{ and }\quad
   X\cap[0\ I_p]\neq\{0\}
 \end{equation}
is an irreducible rational variety of dimension $mp-2$.

Let ${\mathcal X}$ be the set of all $p\times(m+p)$-matrices $X$ whose
entries $x_{i,j}$ satisfy:
$$
\begin{array}{rcl}
x_{i,j}=1&\mbox{ if }&\mbox j=i<p\mbox{ \ or \ } (i,j)=(p,p+1)\\
x_{i,j}\ =\ 0&\mbox{ if }& 
\left\{\begin{array}{ccl}
i=1&\mbox{ and }& j\geq m\\
1<i<p&\mbox{ and }& j<i\mbox{ or } j> i+m\\
i=p&\mbox{ and }& j\leq p\end{array}\right.
\end{array}
$$
The remaining $mp-2$ entries are unconstrained and give 
coordinates for ${\mathcal X}$.
The row space of a matrix $X$ is a $p$-plane $X$ 
satisfying~(\ref{eq:small})
and almost all such $p$-planes arise in this fashion.
Thus ${\mathcal X}$ parameterizes a dense subset of the
subvariety of $p$-planes $X$ satisfying~(\ref{eq:small}).

For example, if $(m,p)=(4,3)$, then ${\mathcal X}$ is the set of all
matrices of the form:
$$
 \left[\begin{array}{ccccccc}
   1   &x_{12}&x_{13}&x_{14}&   0  &   0  &   0  \\
   0   &   1  &x_{23}&x_{24}&x_{25}&x_{26}&   0  \\
   0   &   0  &   0  &   1  &x_{35}&x_{36}&x_{37}
 \end{array}\right].
$$

Since the $1,m$th entry of a matrix in ${\mathcal X}$ vanishes, 
$$
\det\zwei{K(s)}{X}
$$
factors as $s\cdot \psi(s;X)$.
\medskip

\noindent{\bf Conjecture~\ref{conj:shapiroII}$''$. }
{\it
 Let $m,p>1$ be integers.
Then, for almost all choices of non-zero real numbers $s_1,\ldots,s_{mp-2}$,
 the system of equations
\begin{equation}\label{eq:simpler}
\psi(s_1;X)\ =\ \psi(s_2;X)\ =\ \cdots\ =\ 
    \psi(s_{mp-2};X)\ =\ 0
\end{equation}
is zero-dimensional with $d_{m,p}$ solutions, and all of them are real.
}\medskip

The systems of Conjecture~\ref{conj:shapiroII} and the variations given here
are deficient:
They have fewer solutions than standard combinatorial bounds.
For example, the system~(\ref{eq:simpler}) consists of $mp-2$ equations of
degree $p$, thus its B\'ezout number is $p^{mp-2}$.
A better combinatorial bound is 
the normalized volume of the Newton polytope ${\mathcal A}_{m,p}$ of the
polynomial $\psi$~\cite{Kouchnirenko}. 
Table~\ref{table:bounds} compares these combinatorial bounds with $d_{m,p}$, 
for some values of $m,p$.
The volumes of ${\mathcal A}_{m,p}$ were computed using
PHC~\cite{PHC}, a software package for performing general polyhedral
homotopy continuation.
Note the striking difference between the equivalent systems $m,p$ and $p,m$.
\begin{table}[htb]
 \begin{tabular}{|r||c|c|c|c|c|c|c|c|c|c|c|c|}
 \hline
  $m,p$:& 2,2& 3,2& 4,2& 5,2& 6,2&7,2&8,2 & 2,3&3,3 & 4,3&2,4&3,4\\\hline\hline
  $d_{m,p}$:
     &2&5&14&42&132&429&1430&5&42&462&14&462\\\hline
  vol$\,{\mathcal A}_{m,p}$:
     &2&5&18&67&248&919&3426&5&130&3004&42&7156\\\hline
  $p^{mp-2}$: \rule{0pt}{12pt}
     &4&16&64&256&1024&4096&16384&81&2187&59,049&4096&1048576\\\hline
 \end{tabular}\bigskip
 \caption{Combinatorial bounds vs. $d_{m,p}$\label{table:bounds}}
\end{table}

\subsection{Conjecture~\ref{conj:shapiroII} for $m=2, p=3$}

\begin{thm}\label{thm:23}
Conjecture~\ref{conj:shapiroII} holds for $(m,p)=(2,3)$.
\end{thm}

L.~Gonzalez-Vega has also obtained this using resultants and Sturm-Habicht
sequences.
\medskip

\noindent{\bf Proof. }
We will prove the equivalent Conjecture~\ref{conj:shapiroII}$''$.
Let $X:=\{x_{12},x_{23},x_{24},x_{35}\}$ be indeterminates.
Set 
$$
\psi(s;X)\ =\ \ \det \left[\begin{array}{ccccc}
1&s&s^2&s^3&s^4\\
0&1&2s&3s^2&4s^3\\
1&x_{12}&0&0&0\\
0&1&x_{23}&x_{24}&0\\
0&0&0&1&x_{35}\end{array}\right]
$$
We solve the system of polynomials
\begin{equation}\label{eq:param-sys}
   \psi(s;X)\ =\ \psi(t;X)\ =\ 
   \psi(u;X)\ =\ \psi(v;X)\ =\ 0
\end{equation}
by elimination.

The ideal $\langle \psi(s),\psi(t),\psi(u),\psi(v)\rangle$
in the ring ${\mathbb Q}(s,t,u,v)[x_{12},x_{23},x_{24},x_{35}]$ 
has degree $5=d_{2,3}$ and the lexicographic Gr\"obner basis with
$x_{12}<x_{23}<x_{24}<x_{35}$ contains the following univariate polynomial
$g$, which is the universal eliminant for this family of systems:
$$
x_{35}^5  -  4e_1 x_{35}^4  +  (4e_1^2+6e_2) x_{35}^3
- (12e_1e_2+4e_3) x_{35}^2  +  (9e_2^2+8e_1e_3-4e_4) x_{35}
-(12e_2e_3-8e_1e_4)
$$
%
%
Here $e_i$ is the $i$th the elementary symmetric polynomial in $s,t,u,v$.
We show that $g$ has 5 distinct real roots for every choice of distinct
parameters $s,t,u,v$.
The discriminant $\Delta$ of $g$ has degree 20
in the variables $s,t,u,v$ and 711 terms:
$$
\begin{array}{l}
9e_3^4e_2^2e_1^4 - 54e_3^4e_2^3e_1^2 + 81e_3^4e_2^4 - 32e_3^5e_1^5 
+ 204e_3^5e_2e_1^3  - 324e_3^5e_2^2e_1 - 108e_3^6e_1^2 + 324e_3^6e_2 
\\\ \rule{0pt}{13pt}
+ 81e_4^2e_2^4e_1^4 - 486e_4^2e_2^5e_1^2 + 729e_4^2e_2^6 
- 54e_4e_3^2e_2^3e_1^4  + 324e_4e_3^2e_2^4e_1^2 - 486e_4e_3^2e_2^5 
\\\ \rule{0pt}{13pt}
+ 204e_4e_3^3e_2e_1^5 - 1296e_4e_3^3e_2^2e_1^3 
+ 2052e_4e_3^3e_2^3e_1  - 8e_4e_3^4e_1^4 + 738e_4e_3^4e_2e_1^2 
- 2106e_4e_3^4e_2^2 
\\\ \rule{0pt}{13pt}
- 108e_4e_3^5e_1 
- 324e_4^2e_3e_2^2e_1^5 + 2052e_4^2e_3e_2^3e_1^3 - 3240e_4^2e_3e_2^4e_1 
- 108e_4^2e_3^2e_1^6 + 738e_4^2e_3^2e_2e_1^4  
\\\ \rule{0pt}{13pt}
- 2592e_4^2e_3^2e_2^2e_1^2 + 3834e_4^2e_3^2e_2^3 - 368e_4^2e_3^3e_1^3 
+ 1800e_4^2e_3^3e_2e_1 - 27e_4^2e_3^4 + 324e_4^3e_2e_1^6 
\\ \ \rule{0pt}{13pt}
- 2106e_4^3e_2^2e_1^4 + 3834e_4^3e_2^3e_1^2 
- 972e_4^3e_2^4- 108e_4^3e_3e_1^5 + 1800e_4^3e_3e_2e_1^3 
- 5544e_4^3e_3e_2^2e_1 
\\\ \rule{0pt}{13pt}
- 634e_4^3e_3^2e_1^2 + 984e_4^3e_3^2e_2 
- 27e_4^4e_1^4 + 984e_4^4e_2e_1^2 + 432e_4^4e_2^2 
- 352e_4^4e_3e_1 - 64e_4^5.
\end{array}
$$
This vanishes when $g$  has a double root. 
Thus the number of real roots of $g$ is constant on each connected
component (in ${\mathbb R}^4$) of the locus $\Delta\neq 0$.
We show there is only one connected component, and so the number of real
roots of $g$ (and thus the original system)
does not depend upon the choice of real parameters.
Since the roots of $g$  evaluated at $(s,t,u,v)=(1,2,3,4)$
are 
$$
8, \ 8 \pm \sqrt{19}, \ 8 \pm \sqrt{11}
$$
it follows that there are always five real roots of $g$,
and thus the system~(\ref{eq:param-sys}) 
has $d_{2,3}=5$ real solutions whenever $s,t,u,v$ are real and
distinct.

We complete the proof. 
For $w\in{\mathbb Z}_{\geq0}^{10}$, consider the polynomial:
\begin{equation}\label{eq:choice}
s^{w_1}t^{w_2}u^{w_3}v^{w_4}(s-t)^{w_5}(s-u)^{w_6}(s-v)^{w_7}
(t-u)^{w_8}(t-v)^{w_9}(u-v)^{w_{10}}
\end{equation}
Let $A_w$ be the primitive part of the symmetrization of this polynomial.
Thus $A_w$ is a sum of squares, none of which vanish on the locus where
$s,t,u,v$ are distinct.
Then $\Delta$ is 
$$
\begin{array}{c}
\frac{1}{2}(7A_{2220222224}+
3A_{2222402204}+
6A_{4222022222}+
7A_{4220222222}+
2A_{4420022222}\\
\quad +2A_{2222440022}
+2A_{0222443022}
+A_{4420202222}
+2A_{4222420022}
+A_{4220022422}\\
\qquad\hspace{.02in}
+A_{0222442202}
+A_{2202024422}
+6A_{2222420024}
+10A_{4220022242}
+3A_{2222222222}).
\end{array}
$$
Note that the term $7A_{4220222222}$ does not vanish when a single
parameter is zero.
Similarly, the term $3A_{2222402204}$ does not vanish when $s=u$ and $t=v$
(but $u\neq t$).
Thus the locus where $\Delta=0$ has dimension 1 and so its complement is
connected.
\QED

We have a Maple program which performs the computations described and runs
in $\sim$15 seconds on a K6-2-300 processor.

A {\it positive semidefinite} polynomial is a real polynomial that takes only
nonnegative values.  
In  the proof we showed $\Delta$ is positive semidefinite by exhibiting it
as a sum of squares. 
Not all positive semidefinite polynomials are sums of squares of
polynomials.
There exist positive semidefinite polynomials of degree $l$ in $k$
variables which are not sums of squares of polynomials if $\min(k,l)>2$ and 
$(k,l)\neq(3,4)$~\cite{Hilbert_squares}. 
For $\Delta$, $(k,l)=(4,20)$.

The form of the squares we used~(\ref{eq:choice}) for the discriminant
$\Delta$, while motivated by the observation that no two parameters
($0,s,t,u,v,\infty$) should coincide, is justified by the observation that
any real zero of $\Delta$ must also be a zero of all the squares, if
$\Delta$ is a sum of squares.
(See~\cite{CLR_symmetry} for other applications of this idea.)

Each of the polynomials $A_w$ is a sum of squares, the number given by the
orbit of the symmetric group on its index $w$.
Since 6 have trivial stabilizer, 7 are stabilized by a transposition, one by
the dihedral group $D_8$, and one is invariant, there are
$6\cdot 24 + 7\cdot 12 + 3 + 1 = 232$ squares in all.
This is not the best possible.
Choi, Lam, and Reznick~\cite{CLR} show, for degree $l$ homogeneous
polynomials in $k$ variables that are a sum of squares of polynomials, at
most
$$
\Lambda(k,l)\ :=\ \left\lfloor \frac{1}{2} 
\left( \sqrt{1+8\binom{k+l-1}{l}} \ -\ 1\right) \right\rfloor
$$
squares are needed.
Note that $\Lambda(4,20)=59$.

\section{Schubert conditions on a Grassmannian}\label{sec:LR}

\subsection{The Schubert calculus on 
    $\mbox{\it Grass}(p,m+p)$}\label{ssec:LR}
The enumerative problems of Section~\ref{sec:hypersurface} are special cases
of more general problems given by Schubert conditions 
on $\mbox{\it Grass}(p,m+p)$.
A {\it Schubert condition} on $\mbox{\it Grass}(p,m+p)$ is an increasing
sequence of integers 
$$
\alpha\ :\ 1\leq \alpha_1<\alpha_2<\cdots<\alpha_p\leq m+p.
$$
Let $\binom{[m+p]}{p}$ be the set of all such sequences.
A {\it Schubert variety} $\Omega_\alpha\Fdot$ is given by a Schubert
condition $\alpha$ and a complete flag 
$\Kdot$ in ${\mathbb C}^{m+p}$, a sequence of subspaces
$$
\Kdot\ :\ K_1\subset K_2\subset \cdots \subset K_{m+p}={\mathbb C}^{m+p}
$$
where $\dim K_i=i$.
Then the Schubert variety $\Omega_\alpha\Fdot$ is the set of all 
$p$-planes $X$ satisfying
\begin{equation}\label{eq:schubert_condition}
  \dim X\cap K_{m+p+1-\alpha_i}\ \geq \ p+1-i
\end{equation}
for each $i=1,2,\ldots,p$.
This irreducible subvariety of 
$\mbox{\it Grass}(p,m+p)$ has 
codimension $|\alpha|:=\sum_i (\alpha_i-i)$.

A sequence
$\alpha^\DOT=\alpha^1,\ldots,\alpha^n$  with $\alpha^j\in\binom{[m+p]}{p}$
with $\sum_j|\alpha^j|=mp$
is {\it Schubert data} for $\mbox{\it Grass}(p,m+p)$.
Given Schubert data $\alpha^\DOT$ 
and flags $K^1_\DOT,\ldots,K^n_\DOT$ in general position, 
there are finitely many (complex)
$p$-planes $X$ which lie in the intersection of the Schubert 
varieties $\Omega_{\alpha^j}K^j_\DOT$ for $j=1,\ldots,n$. 
The classical Schubert calculus~\cite{KlLa} gives 
the following recipe for computing this number  $d=d(m,p;\alpha^\DOT)$.
Let $h_1,\ldots,h_m$ be indeterminates with $degree(h_i) = i$.
For each integer sequence $\beta_1<\beta_2<\cdots<\beta_r$
define the following polynomial:
\begin{equation}\label{schur}
S_\beta\ :=\ \det ( h_{\beta_i-j})_{1\leq i,j\leq r}.
\end{equation}
Here $h_0 := 1$ and $h_i:=0$ if $i<0$ or $i>m$.
Let ${\mathcal I}$ be the ideal in 
${\bf Q}[h_1,\ldots,h_m]$ generated by those 
$S_\beta$ with  $r=p+1$, $1<\beta_1$, and $\beta_{p+1}\leq m+p$.
The quotient ring 
${\mathcal A}_{m,p}:= {\bf Q}[h_1,\ldots,h_m]/{\mathcal I}$
is isomorphic to the cohomology ring of $\mbox{\it Grass}(p,m+p)$. 
It is Artinian with one-dimensional socle in degree $mp$.
In the socle we have the relation
$$
 d \cdot (h_m)^p \, -  \,
S_{\alpha^1} S_{\alpha^2} \cdots S_{\alpha^n}  \quad \in \quad {\mathcal I}. 
$$
We can compute the number $d$ by normal form reduction modulo 
any Gr\"obner basis for ${\mathcal I}$. 

If $\gamma$ is a rational normal curve, then the flag of subspaces
osculating $\gamma$ at a point is the {\it osculating flag} to $\gamma$ at
that point.

\begin{conj}\label{conj:shapiro-grass}
Let $m,p>1$ and $\alpha^\DOT$ be Schubert data for 
$\mbox{\it Grass}(p,m+p)$.
For almost all choices of 
flags  $K^1_\DOT,\ldots,K^n_\DOT$ osculating a fixed rational
normal curve at real points, there
are exactly $d(m,p,\alpha^\DOT)$ $p$-planes $X$ in the intersection
of Schubert varieties
$$
\Omega_{\alpha^1}K^1_\DOT \cap \Omega_{\alpha^2}K^2_\DOT\cap
 \cdots\cap \Omega_{\alpha^n}K^n_\DOT,
$$
and each of these $p$-planes is real.
\end{conj}

As with Conjecture~\ref{conj:shapiroII}, the intersection is
zero-dimensional if the points of osculation are distinct~\cite{EH83},
and there are no multiplicities for the
important class of Pieri Schubert data, (described below) which includes the
case of Conjecture~\ref{conj:shapiroII}. 

If $\alpha_i=1+\alpha_{i-1}$, then 
condition~(\ref{eq:schubert_condition}) for $i-1$
implies~(\ref{eq:schubert_condition}) for $i$. 
Thus only those conditions~(\ref{eq:schubert_condition})
with $\alpha_i-\alpha_{i-1}>1$ (or $\alpha_1>1$) are {\it essential},
and so only the subspaces $K_{m+p+1-\alpha_i}$ corresponding to essential
conditions need be specified in a flag.
If $\alpha:=(1,2,\ldots,p-1,p+1)$, then only the last condition is
essential, thus the Schubert variety $\Omega_\alpha \Kdot$ consists of
those $X$ with $\dim X\cap K_m\geq 1$.
This shows Conjecture~\ref{conj:shapiroII} is a special case of 
Conjecture~\ref{conj:shapiro-grass}.

\subsection{Systems of polynomials}\label{sec:loc-coords}
A complete flag $\Kdot$ is represented by a nonsingular matrix
also written $\Kdot$:
The $i$-plane $K_i$ is the row space of $K_i$, the first $i$ rows of
$\Kdot$.
The condition that $\dim X\cap K_{m+p+1-\alpha_i}\geq p+1-i$ is given
by 
$$
\mbox{$(m+p+1+i-\alpha_i)$-minors of }\ 
\zwei{K_{m+p+1-\alpha_i}}{X}\ =\ 0.
$$

The flag $\Kdot(s)$ osculating the rational normal curve
$\gamma$ with the parameterization~(\ref{eq:afine_rat}) at $\gamma(s)$ is
represented by the  $(m+p)\times(m+p)$-matrix whose $i,j$th entry is
$\binom{j-i}{i-1}s^{j-i}$.
\medskip

\noindent{\bf Conjecture~\ref{conj:shapiro-grass}$'$. }
{\it 
Let $m,p>1$ and $\alpha^\DOT$ be Schubert data for 
$\mbox{\it Grass}(p,m+p)$.
For almost all $n$-tuples of distinct real numbers $s_1,\ldots,s_n$,
the system of polynomials 
$$
(m+p+1+i-\alpha_i^j)\mbox{-minors of }\ 
\zwei{K_{m+p+1-\alpha_i^j}(s_j)}{I_p\ \ X}\ =\ 0
$$
for $i=1,\ldots,p$ and $j=1,\ldots,n$ has 
$d(m,p,\alpha^\DOT)$ solutions, and each is real.
}\medskip

For any Schubert conditions 
$\alpha,\beta$ with $\alpha_i+\beta_{p+1-i}\leq m+p$
for $i=1,\ldots,p$, let ${\mathcal X}_{\alpha,\beta}$ be the collection of
all $p\times(m+p)$-matrices $X$ whose entries $x_{ij}$ satisfy
\begin{equation}\label{eq:loc_coords}
 \begin{array}{rcl}
  x_{i,\alpha_i}\ =\ 1&\quad&\mbox{for }i=1,\ldots,p\\
  x_{i,j}\ =\ 0&&\mbox{if }j<\alpha_i\ \mbox{ or }\ j>m+p+1-\beta_{p+1-i}
 \end{array}
\end{equation}

If $X\in{\mathcal X}_{\alpha,\beta}$, then the row space of $X$ is  a
$p$-plane in the intersection 
$\Omega_\alpha\Kdot(\infty)\cap \Omega_\beta\Kdot(0)$.
In this way, ${\mathcal X}_{\alpha,\beta}$ parameterizes a Zariski open
subset of the set of all such $p$-planes.
This parameterization can be used to obtain a system of equations simpler
than, but equivalent to, the system of
Conjecture~\ref{conj:shapiro-grass}$'$.

The map ${\mathcal X}_{\alpha,\beta}\rightarrow\mbox{Grass}(p,m+p)$ is not
injective.
For example, ${\mathcal X}_{123,134}$ consists of all $3\times 7$-matrices
of the form:
$$
\left[
\begin{array}{ccccccc}
1&x_{12}&x_{13}&x_{14}&0&0&0\\
0&1&x_{23}&x_{24}&x_{25}&0&0\\
0&0&1&x_{34}&x_{35}&x_{36}&x_{37}\end{array}\right]
$$
Let $r_1,r_2,r_3$ be the rows of such a matrix.
If $x_{36}=x_{37}=0$, then for each $a\in{\mathbb C}$, the matrix with rows  
$r_1,r_2+ar_3,r_3$ is in  ${\mathcal X}_{123,134}$, and these all have the
same row space.
Similarly, if $x_{25}=0$, then the same is true of the matrices with rows
$r_1+ar_2,r_2,r_3$.

Let ${\mathcal X}^\circ_{\alpha,\beta}\subset{\mathcal X}_{\alpha,\beta}$
be the set of those matrices whose entries further satisfy
\begin{equation}\label{eq:non-zero}
 \mbox{
  \begin{minipage}{4.5in}
  For each $i=2,\ldots,p$, at least one $x_{ij}\neq 0$, for $j$ satisfying 
  $$
    \beta_{p+1-i}\leq  m+p+1-j < \beta_{p+2-i}.
  $$
  \end{minipage}
 }
\end{equation}
The map 
${\mathcal X}^\circ_{\alpha,\beta}\rightarrow\mbox{Grass}(p,m+p)$ is
injective. 
For ${\mathcal H}_{123,134}$ this condition is that $x_{25}\neq 0$ and
$(x_{36},x_{37})\neq(0,0)$.

\subsection{Pieri Schubert conditions}

If $\alpha\in\binom{[m+p]}{p}$ has $\alpha_{p-1}=p-1$ and $\alpha_p=p+a$,
then the Schubert variety $\Omega_\alpha\Kdot$ is
$$
\{X \mid X\cap K_{m+1-a}\neq \{0\}\}.
$$
We call such a Schubert condition a {\it Pieri condition} and
denote it by $J_a$.
{\it Pieri Schubert data} are Schubert data $\alpha^1,\ldots,\alpha^n$
were at most 2 of the conditions $\alpha^i$ are {\it not} Pieri
conditions. 
These include the Schubert data of Conjecture~\ref{conj:shapiroII}.

\begin{prop}[Theorem 9.1 in~\cite{EH83}]\label{no-mult}
If $\alpha^\DOT$ are Pieri Schubert data and the flags 
$K_\DOT^1,\ldots,K_\DOT^n$ osculate a rational normal curve at general
points, then the intersection of Schubert varieties
$$
\Omega_{\alpha^1}K^1_\DOT \cap \Omega_{\alpha^2}K^2_\DOT\cap
 \cdots\cap \Omega_{\alpha^n}K^n_\DOT,
$$
is transverse.
In particular, there are no multiplicities.
\end{prop}

Here is the main theorem of this section.

\begin{thm}\label{thm:implies}
Let $a,b>1$ and suppose that Conjecture~\ref{conj:shapiroII} holds for this
$(m,p)=(a,b)$.
Then Conjecture~\ref{conj:shapiro-grass} holds for any Pieri Schubert
data for
$\mbox{\it Grass}(p,m+p)$ with $(p,m)\leq (a,b)$ or $(p,m)\leq (b,a)$,
in each coordinate. 
\end{thm}

\begin{rem}\label{rem:stronger}
If the conclusion of Proposition~\ref{no-mult} held for all Schubert data,
then the proof we give of Theorem~\ref{thm:implies} would imply its
conclusion for all Schubert data as well.
\end{rem}

Pieri conditions are special because of Pieri's formula.
For $\alpha,\beta\in\binom{[m+p]}{p}$ and $a>0$, we write
$\alpha<_a\beta$ if $|\alpha|+a=|\beta|$ and 
$$
\alpha_1\leq\beta_1<\alpha_2\leq\beta_2<\cdots<\alpha_p\leq \beta_p.
$$

\begin{prop}[Pieri's Formula]\label{prop:pieri}
Let $J_a:=1<2<\cdots<p{-}1<p{+}a\in\binom{[m+p]}{p}$.
\begin{enumerate}
\item[{\rm (i)}]
In the cohomology ring ring ${\mathcal A}_{m,p}$ of 
$\mbox{\it Grass}(p,m+p)$, $S_{J_a}=h_a$ and
$$
S_\alpha \cdot S_{J_a}\ =\ \sum_{\alpha<_a\beta} S_\beta.
$$

\item[{\rm (ii)}]
If $\Kdot(s)$ and $\Kdot(t)$ are flags osculating a rational normal curve at 
points $s$ and $t$, then 
$$
\lim_{s\rightarrow t}\left(
\Omega_\alpha\Kdot(t) \cap \Omega_{J_a}\Kdot(s)\right)\ =\ 
\sum_{\alpha<_a\beta} \Omega_\beta\Kdot(t).
$$
Here, the limit is taken as cycles.
By this we mean that the sum is the fundamental cycle of the
limit of the schemes $\Omega_\alpha\Kdot(t) \cap \Omega_{J_a}\Kdot(s)$
as $s$ approaches $t$ along the rational normal curve.

\item[{\rm (iii)}]
Suppose $\alpha^\DOT=\alpha^1,J_a,\alpha^2,\ldots,\alpha^n$ are Schubert
data.
Then
$$
d(m,p;\alpha^\DOT)\ =\ \sum_{\alpha<_a\beta} 
d(m,p;\beta,\alpha^2,\ldots,\alpha^n).
$$
\end{enumerate}
\end{prop}

Statement (i) is the usual statement of Pieri's
formula~\cite{Fulton_tableaux,Hodge_Pedoe},
Statement (ii) is Theorem 8.1 of~\cite{EH83}, and Statement (iii)
is a direct consequence of (i).
\medskip

\noindent{\bf Proof.}
We deduce Theorem~\ref{thm:implies} from Lemma~\ref{lem:ordering} below,
which shows some simple dependencies between
Conjecture~\ref{conj:shapiro-grass} for different collections of Schubert
data.
Definition~(\ref{eq:schubert_condition}) implies that 
$\Omega_\beta\Kdot\subset \Omega_\alpha\Kdot$ if and only if 
$\alpha\leq \beta$ coordinatewise. 
In fact, $\Omega_\beta\Kdot\cap
\Omega_\alpha\Kdot=\Omega_{\beta\vee\alpha}\Kdot$, where $\beta\vee\alpha$
is the coordinatewise maximum of $\alpha$ and $\beta$.
We make some definitions needed for the statement of Lemma~\ref{lem:ordering}.

\begin{defn}\label{def:triple}
 Let $m,p>1$ be integers. 
\begin{enumerate}
\item   
     For $\alpha\in\binom{[m+p]}{p}$ define
     $\alpha^\perp\in\binom{[m+p]}{m}$ to 
     be the increasing sequence obtained from the numbers
     $\{1,2,\ldots,m+p\}\setminus\{\alpha_1,\ldots,\alpha_p\}$.
     Given Schubert data $\alpha^\DOT$ for $\mbox{\it Grass}(p,m+p)$, set
     $\alpha^{\DOT\perp}$ to be $(\alpha^1)^\perp,\ldots,(\alpha^n)^\perp$.
\item
     Suppose $p>2$.
     For $\alpha\in\binom{[m+p-1]}{p-1}$ define
     $\alpha^+\in\binom{[m+p]}{p}$ to 
     be $1<1+\alpha_1<\cdots<1+\alpha_{p-1}$. 
     Given Schubert data $\alpha^\DOT$ for $\mbox{\it Grass}(p,m+p)$, set
     $\alpha^{\DOT+}$ to be $(\alpha^1)^+,\ldots,(\alpha^n)^+$.
\item
     Let $\preceq$ be the partial order on Pieri Schubert data where we
     say that  $\beta^\DOT$ covers
     $\alpha^\DOT=\alpha^1,\ldots,\alpha^n$ if one of the following holds
     $$
     \begin{array}{cc}
     \beta^\DOT\ =\ \beta,\alpha^3,\ldots,\alpha^n\ \mbox{with}\
     \alpha^2=J_a\ \mbox{and}\ \alpha^1<_a\beta\\
     \beta^\DOT\ =\ \alpha^1,\ldots,\alpha^{n-2},\beta\ \mbox{with}\
     \alpha^{n-1}=J_a\ \mbox{and}\ \alpha^n<_a\beta.
     \end{array}
     $$
\end{enumerate}
\end{defn}

\begin{lem}\label{lem:ordering}
Let $m,p>1$ be integers.
\begin{enumerate}
\item[(i)]
     If $\alpha^\DOT$ is Schubert data for $\mbox{\it Grass}(p,m+p)$, then
     $\alpha^{\DOT\perp}$ is Schubert data for $\mbox{\it Grass}(m,m+p)$.
     Moreover, Conjecture~\ref{conj:shapiro-grass} holds for
     $m,p,\alpha^\DOT$ if and only if it holds for $p,m,\alpha^{\DOT\perp}$.
\item[(ii)]
     Suppose $p>2$ and let $J_m:=1<2<\cdots<p-1<p+m$.
     If $\alpha^\DOT$ is Schubert data for $\mbox{\it Grass}(p-1,m+p-1)$,
     then $\beta^\DOT:=\alpha^{\DOT+},J_m$ is Schubert data for
     $\mbox{\it Grass}(m,m+p)$. 
     Moreover, Conjecture~\ref{conj:shapiro-grass} holds for
     $m,p-1,\alpha^\DOT$ if and only if it holds for $m,p,\beta^{\DOT}$.
\item[(iii)]
     Let $\alpha^\DOT,\beta^\DOT$ be Pieri Schubert data for 
     $\mbox{\it Grass}(p,m+p)$ with $\alpha^\DOT\preceq\beta^\DOT$.
     If Conjecture~\ref{conj:shapiro-grass} holds for
     $\alpha^\DOT$ for $\mbox{Grass}(p,m+p)$, then it holds for 
     $\beta^\DOT$.
\end{enumerate}
\end{lem}

\noindent{\bf Proof of Theorem~\ref{thm:implies}. }
First note that Conjecture~\ref{conj:shapiro-grass} holds for Schubert data
$\alpha^\DOT$ for $\mbox{\it Grass}(p,m+p)$ if and only if it holds for any
rearrangement of the data $\alpha^\DOT$.
Suppose Conjecture~\ref{conj:shapiroII} holds for $\mbox{Grass}(b,a+b)$.
Let $\alpha^\DOT$ be Pieri Schubert data for $\mbox{\it Grass}(p,m+p)$ where
$(m,p)\leq(a,b)$ or $(m,p)\leq(b,a)$ coordinatewise.
Since $J_1^\perp=J_1$, Conjecture~\ref{conj:shapiroII} holds also for
$(m,p)=(b,a)$, by Lemma~\ref{lem:ordering}(i). 
Thus we may assume that $(m,p)\leq(a,b)$.
By Lemma~\ref{lem:ordering}(ii), there exist Pieri Schubert data $\beta^\DOT$
for $\mbox{Grass}(b,a+b)$ such that Conjecture~\ref{conj:shapiro-grass}
holds for $\alpha^\DOT$ if and only if it holds for $\beta^\DOT$.
Finally, Theorem~\ref{thm:implies} follows from (iii) by noting that the
Schubert data of Conjecture~\ref{conj:shapiroII}, namely 
$\alpha^1=\cdots=\alpha^{ab}=J_1$, is minimal among all
Pieri Schubert data for $\mbox{\it Grass}(b,a+b)$.
\QED

\noindent{\bf Proof of Lemma~\ref{lem:ordering}. }
For (i), fix a real inner inner product on ${\mathbb C}^{m+p}$.
Then the map $X\mapsto X^\perp$ gives an isomorphism between
$\mbox{\it Grass}(p,m+p)$ and 
$\mbox{\it Grass}(m,m+p)$.
Given a flag $\Kdot$ and an increasing sequence $\alpha$, let 
$\Kpdot$
be the flag of annihilators of the subspaces of $\Kdot$.
Then we have
$$
X\in \Omega_\alpha\Kdot\ \Longleftrightarrow\  
X^\perp\in \Omega_{\alpha^\perp}\Kpdot.
$$
Furthermore, if $\Kdot(s)$ is the flag of subspaces osculating a rational
normal curve $\gamma$ at a point $\gamma(s)$, then $(K_{m+p-1}(s))^\perp$ is
a rational normal curve with $\Kpdot(s)$ its osculating flag.
Thus Conjecture~\ref{conj:shapiro-grass} for Schubert data $\alpha^\DOT$ for
$\mbox{Grass}(p,m+p)$ is equivalent to 
Conjecture~\ref{conj:shapiro-grass} for Schubert data $\alpha^{\DOT\perp}$
for $\mbox{Grass}(m,m+p)$.
\smallskip

For (ii), let $\gamma$ be the rational curve~(\ref{eq:afine_rat}) with
$\Kdot(s)$ as before. 
Then $X\in \Omega_{J_m}\Kdot(\infty)$ if and only if 
$\langle\gamma(\infty)\rangle=K_1(\infty)\subset X$.
Consider the projection 
$\pi:{\mathbb C}^{m+p}\twoheadrightarrow{\mathbb C}^{m+p-1}$ from the last
coordinate $\gamma(\infty)$.
If $X\in \Omega_{J_m}\Kdot(\infty)$, then $X':=\pi X$ is a $(p{-}1)$-plane.
This induces an isomorphism 
$\pi:\Omega_{J_m}\Kdot(\infty)\stackrel{\sim}{\longrightarrow}
\mbox{\it Grass}(p{-}1,m{+}p{-}1)$.
The inverse map is given by $X'\mapsto K_1(\infty)+X'$.

The projection $\pi\circ\gamma$ is the standard rational normal curve
$\gamma'$ in  ${\mathbb C}^{m+p-1}$.
Similarly, the flag $\Kdot'(s)$ osculating $\gamma'$ at $\gamma'(s)$ is
$\pi \Kdot(s)$.
Note that if $L$ is a linear subspace of ${\mathbb C}^{m+p}$ 
with $\gamma(\infty)\not\in L$, then $\dim X\cap L=\dim \pi X\cap \pi L$.
In particular, if $X\in \Omega_{J_m}\Kdot(\infty)$, $s\neq \infty$,  and 
$\alpha\in\binom{[m+p-1]}{p-1}$, then
$\dim  X'\cap K'_{(m+p-1)+1-\alpha_i}\geq (p-1)+1-i$ if and only if 
$\dim X\cap K_{m+p+1-(1+\alpha_i)}\geq p+1-(i+1)$.
Thus we have
\begin{equation}\label{eq:equiv-proj}
  X\in \Omega_{J_m}\Kdot(\infty)\cap \Omega_{\alpha^+}\Kdot(s)
  \ \Longleftrightarrow\ X'\in \Omega_{\alpha}\Kdot'(s).
\end{equation}
In fact, this induces an isomorphism of schemes.

This gives a strong equivalence between enumerative problems:
If $\alpha^1,\ldots,\alpha^n$ are in $\binom{[m+p-1]}{p-1}$ and
$s_1,\ldots,s_n$ any complex numbers, then 
the map $\pi$ induces an isomorphism between the schemes
$$
\Omega_{J_m}\Kdot(\infty)\cap\bigcap_{i=1}^n \Omega_{(\alpha^i)^+}\Kdot(s_i)
\quad\mbox{ and }\quad
\bigcap_{i=1}^n \Omega_{\alpha^i}\Kdot'(s_i).
$$
Part (ii) follows by noting that any real reparameterization of the rational
normal curve $\gamma$ induces an 
isomorphism of polynomial systems, thus preserves real solutions.
Hence given $s_0,s_1,\ldots,s_n\in{\mathbb P}^1_{\mathbb R}$, there is 
an equivalent system with $s_0=\infty$. 
\medskip

It suffices to prove (iii) when $\beta^\DOT$ covers
$\alpha^\DOT$ in the partial order $\preceq$ defined on Pieri Schubert data.
Suppose Conjecture~\ref{conj:shapiro-grass} fails for $\beta^\DOT$
and $\beta^\DOT$ covers $\alpha^\DOT$ with $\alpha^1<_a\beta$ and
$\alpha^2=J_a$ as in Definition~\ref{def:triple} (iii).
Then there exist distinct real numbers $s_1,s_3,\ldots,s_n$
such that 
\begin{equation}\label{eq:intersec}
 \Omega_{\beta}\Kdot(s_1) \cap 
\Omega_{\alpha^3}\Kdot(s_3) \cap \cdots \cap \Omega_{\alpha^n}\Kdot(s_n)
\end{equation}
is transverse with some complex $p$-planes in the intersection.
We may assume without any loss that $s_1=0$.
Then there is an open subset ${\mathcal O}$ of the set of
$(n-1)$-tuples of real numbers $s_3,\ldots,s_n$ such
that~(\ref{eq:intersec}) is transverse and contains a complex $p$-plane
$X$.

By the dimensional transversality results of~\cite{EH83},   
we may assume further that for $\beta'\in\binom{[m+p]}{p}$ and
$(s_2,\ldots,s_n)\in {\mathcal O}$, the
intersection
$$
\Omega_{\beta'}\Kdot(0)\cap \bigcap_{i=3}^{n}\Omega_{\alpha^i}\Kdot(s_i) 
$$
has the expected dimension and is transverse if 0-dimensional.
This is empty if $|\beta'|>a+|\alpha^1|$, for dimension
reasons. 
Thus 
$$
\left(\sum_{\alpha<_a\beta'}\Omega_{\beta'}\Kdot(0)\right)
\cap \bigcap_{i=3}^{n}\Omega_{\alpha^i}\Kdot(s_i) 
$$
is transverse for $(s_3,\ldots,s_n)\in {\mathcal O}$.

Fix $(s_3,\ldots,s_n)\in {\mathcal O}$.
By Proposition~\ref{prop:pieri}(i), 
there is an $\epsilon>0$ 
such that for $|t|<\epsilon$
$$
\Omega_\alpha\Kdot(0) \cap \Omega_{J_a}\Kdot(t) 
\cap \bigcap_{i=3}^{n}\Omega_{\alpha^i}\Kdot(s_i) 
$$
is transverse.
Here, when $t=0$, replace $\Omega_\alpha\Kdot(0) \cap \Omega_{J_a}\Kdot(t)$
by $\sum_{\alpha<_a\beta}\Omega_{\beta'}\Kdot(0)$.
Since at $t=0$ not all points in the intersection are real, the same holds
for $0<t<\epsilon$.
But then Conjecture~\ref{conj:shapiro-grass} fails for 
the Schubert data $\alpha^\DOT$
and completes the proof of Lemma~\ref{lem:ordering}.
\QED

\subsection{An infinite family}

We show that Conjecture~\ref{conj:shapiro-grass} holds for an infinite
family of nontrivial Schubert data.

\begin{thm}\label{infinite_family}
Conjecture~\ref{conj:shapiro-grass} holds for any $m$ with  
$p=2$ and Pieri Schubert data where one condition is $J_{m-1}$.
\end{thm}

\noindent{\bf Proof. }
By Lemma~\ref{lem:ordering}(iii), it suffices to show this for
$\alpha^1=\cdots=\alpha^{m+1}=J_1$ and $\alpha^{m+1}=J_{m-1}$.
Geometrically, we are looking for the 2-planes which meet a 2-plane
and $m+1$ general $m$-planes nontrivially. 
We first show there are $m$ such 2-planes.
Let $L=K_2(\infty)=[0\ I_2]$ and $M=K_m(0)=[I_m\ 0]$, and let
$N_i=K_m(s_i)$, 
where $s_1,\ldots,s_m$ are distinct nonzero real numbers.
For each one-dimensional subspace $\lambda$ of $L$ and each 
$1\leq i\leq m$, the composition
$$
M\ \hookrightarrow\  L\oplus M \simeq {\mathbb C}^{m+2}\ 
\twoheadrightarrow\ 
L\oplus M/(\lambda+N_i)\ \simeq {\mathbb C}
$$
defines a linear form $\psi_{i,\lambda}$ on $M$.
Each one-dimensional subspace $\mu$ of its kernel gives a 2-plane
$\lambda\oplus\mu$ containing $\lambda$ and meeting both $M$ and $N_i$
nontrivially. 

Thus if $X$ is a 2-plane meeting $L,M$, and each $N_i$ nontrivially, 
then $H\cap L=\lambda$ and $H\cap M=\mu$ are lines with $\mu$ in the kernel
of each form $\psi_{i,\lambda}$.
Hence the forms are dependent.
Similarly, if $\lambda$ is a line in $L$ such that the forms 
$\psi_{i,\lambda}$ are dependent, then any line $\mu$ they
collectively annihilate gives a 2-plane $\lambda\oplus\mu$ meeting
$L,M$, and each $N_i$ nontrivially.
It follows that the number of such 2-planes is the degree of the determinant
of the forms $\psi_{i,\lambda}$, a polynomial in 
$\lambda\in{\mathbb P}(L)\simeq{\mathbb P}^1$.
Since each form $\psi_{i,\lambda}$ is a linear function of
$\lambda$, the determinant has degree $m$, so there are $m$ 2-planes $X$
meeting $L,M$, and each $N_i$ nontrivially.

We compute this determinant and show it has only real roots.
Let $\lambda=\lambda(x)$ be the span of the vector
$$
(0,\ldots,0,1,(m+1)x).
$$
Let the rational normal curve $\gamma$ have the parameterization
$$
\gamma\ :\ s\ \longmapsto\ (1,-s,s^2,\ldots,(-1)^{m+1}s^{m+1}).
$$
Then $K_m(s)$, the osculating $m$-plane to $\gamma$ at $\gamma(s)$, 
is the kernel of the matrix:
\begin{equation}\label{vf}
  \left[\begin{array}{cccccccc}
  s^m&ms^{m-1}&\ldots&  {m\choose{j}}s^{m-j}  &\ldots&    ms         &1&0\\
  0  &   s^m  &\ldots&{m\choose{j-1}}s^{m-j+1}&\ldots&\binom{m}{2}s^2&ms&1
  \rule{0pt}{14pt}
  \end{array}
  \right].
\end{equation}
If $R_j(s)$ is the linear form given by the $j$th row of this 
matrix, then 
$$
((m+1)x+ms_i)R_1(s_i)\ -\ R_2(s_i)
$$ 
vanishes on $\lambda(x)$ and its restriction to $M$ gives the form
$\psi_{i,\lambda(x)}$.
This restriction is represented by the vector $\Lambda(s_i,x)$ 
whose $j$th coordinate for $j=0,\ldots,m-1$ is:
$$
{\textstyle {m+1\choose j}}
\left((m{-}j{+}1)xs_i^{m-j}+ (m{-}j)s_i^{m-j+1}\right).
$$
We seek the determinant of the following matrix:
$$
\left[\begin{array}{c}
\Lambda(s_1,x)\\ \vdots\\\Lambda(s_m,x)\end{array}\right].
$$
This factors as $A\cdot B$, where $A$ is the
bidiagonal $m\times (m+1)$-matrix
$$
\left[
\begin{array}{cccccccccc}
m&(m+1)x&0\\
0&m^2-1&m(m+1)x&0\\
&\ \makebox[.1in][l]{\qquad$\ddots$}&
           \makebox[.1in][l]{\ \qquad$\ddots$}&\\
&0&{m+1\choose j}(m{-}j)& {m+1\choose j}(m{-}j{+}1)x&&0\\
&&\makebox[.1in][l]{\qquad\quad$\ddots$}&
\makebox[.1in][l]{\qquad\quad\ $\ddots$}&\\
&&0&{m+1\choose 2}\quad&&{m+1\choose 2}2x\end{array}\right]
$$
and $B$ is the $(m+1)\times m$-matrix whose $i,j$th entry is
$s_j^{m+2-i}$.
Numbering the rows of $A$ and the columns of $B$ from 0 to $m$, we see
that
$$
\det(A(x)\cdot B)\ =\ \
        \sum_{i=0}^m (-1)^i \det A_i(x) \det B_i,
$$
where $A_i$ is the matrix $A$ with its $i$th column removed and $B_i$
is the matrix $B$ with its $i$th row removed.
We find that 
\begin{eqnarray*}
\det(A_i)&=& m!(m+1-i)x^{m-i}\prod_{j=1}^m{\textstyle \binom{m+1}{j}}\\
\det B_i&=& e_i(s_1,\ldots,s_m)s_1s_2\cdots s_m\cdot\prod_{j<k}(s_j-s_k)
\end{eqnarray*}
and so $\det(A\cdot B)$ is
$$
m! \prod_{j<k}(s_j-s_k)\prod_{j=1}^m s_j {\textstyle \binom{m+1}{j}}
\cdot \left(\sum_{i=0}^m (-1)^i(m-i+1)x^{m-i}e_i(s_1,\ldots,s_m)\right).
$$

Thus the coordinate $x$ of the line $\lambda$ satisfies the polynomial
$$
P_m(s_1,\ldots,s_m;x)\ :=\ 
\sum_{i=0}^m (-1)^i(m-i+1)x^{m-i}e_i(s_1,\ldots,s_m).
$$
Since we have 
$e_i(s_1,\ldots,s_m)=e_i(s_1,\ldots,s_{m-1})+s_me_{i-1}(s_1,\ldots,s_{m-1})$,
we see that 
$$
P_m(s_1,\ldots,s_m;x)= (x-s_m)P_{m-1}(s_1,\ldots,s_{m-1};x) +
 x \prod_{i=1}^{m-1}(x-s_i).
$$

To complete the proof, we use induction to show that
$$
(*)\qquad\qquad \mbox{
\begin{minipage}{4.5in}
If $0<s_1<\cdots<s_m$, then the roots $r_1,\ldots,r_m$
of $P_m$ satisfy:
$$
0<r_1<s_1<r_2<s_2<\cdots<r_m<s_m.
$$
\end{minipage}
}\qquad
$$
This suffices, if we can assume $0<s_1<\cdots<s_m$.
But we may assume this:
Given a set of distinct real numbers  
$s_1,\ldots,s_m,s_{m+1},s_{m+2}$, we may assume $s_{m+2}=\infty$ and 
$s_{m+1}<s_1<\cdots<s_m$ and then apply the automorphism 
$s \mapsto s-s_{m+1}$ of ${\mathbb P}^1({\mathbb R})$ which fixes
$\infty=s_{m+2}$.

The case $m=1$ of $(*)$ holds as $P_1(s_1;x)=2x-s_1$.
Suppose $P_{m-1}$ satisfies $(*)$.
Then the roots of $(x-s_m)P_{m-1}$ are $r_1<r_2<\cdots<r_{m-1}<s_m$
and those of $x\prod_{i=1}^{m-1}(x-s_i)$ are
$0<s_1<\cdots<s_{m-1}$.
Moreover the leading coefficients of both polynomials are positive.
The result follows by the Intermediate Value Theorem:
If $P(x)$ and $Q(x)$ are polynomials of degree $n$ with positive leading
coefficients and real interlaced roots $p_i$ of $P$ and $q_i$ of $Q$
$$
p_1<q_1<p_2<q_2<\cdots<p_n<q_n,
$$
then $P(x)+Q(x)$ has real roots $r_i$ satisfying
$p_i<r_i<q_i$, for $i=1,\ldots, n$.
\QED

\subsection{Computational evidence}
We have proven Conjecture~\ref{conj:shapiro-grass} in a number of cases
besides those of Theorem~\ref{infinite_family}.
We also have done many computations along the lines of those in
Section~\ref{sec:computation}.
To describe these, we use the following compact notation.
If a Schubert condition $\alpha$ is repeated $k$ times in some Schubert
data, we abbreviate that by $\alpha^k$.
Thus, the conditions of Conjecture~\ref{conj:shapiroII} are written
as $J_1^{mp}$.

\begin{thm}\label{thm:grass-proof}
Conjecture~\ref{conj:shapiro-grass} holds for the following 
Schubert data.
\begin{enumerate}

\item[(i)] $(m,p)=(4,2)$, $\alpha^\DOT=J_2^4$.
	Here, $d(4,2;J_2^4)=3$.	

\item[(ii)]  $(m,p)=(3,3)$, $\alpha^\DOT=J_2^4,J_1$.
	Here, $d(3,3;J_2^4,J_1) = 3$.

\item[(iii)]  $(m,p)=(3,3)$, $\alpha^\DOT=(135)^2,J_1^3$.
	Here, $d(3,3;(135)^2,J_1^3) = 6$.

\item[(iv)] $(m,p)=(4,3)$, $\alpha^\DOT=135^4$.  
	Here, $d(4,3;135^4)=8$.

\end{enumerate}
\end{thm}

\noindent{\bf Proof. }
We consider a polynomial system with parameters, give a universal eliminant,
and show the eliminant has only real roots for distinct values of the
parameters.
We work in the local parameterization ${\mathcal X}_{\alpha^1,\alpha^2}$
of Section~\ref{sec:loc-coords}.

\noindent{(i) }
Let $(m,p)=(4,2)$ and  $\alpha^\DOT=J_2^4$.
The equations are
$$
\mbox{maximal minors }
\left[\begin{array}{cccccc}
1&s&s^2&s^3&s^4&s^5\\
0&1&2s&3s^2&4s^3&5s^4\\
0&0&1&3s&6s^2&10s^3\\
1&x_{12}&x_{13}&0&0&0\\
0&0&0&1&x_{25}&x_{26}\end{array}\right]
\ =\ 0
$$
and the same equations with $t$ replacing $s$.
The ideal of these polynomials contains the following univariate polynomial
$g$ of degree $3=d(4,2,J_2^4)$.
$$
25x_{12}^3-25x_{12}^2(s+t) +x_{12}(19st+6s^2+6t^2)-3(s^2t+st^2)
$$
whose discriminant has primitive part
$$
9(s-t)^6 + 23s^2t^2(s-t)^2 + 9(s^6+t^6).
$$
Since $g(x_{12};1,2)$ has roots
$$ 1,\ 1\pm\frac{1}{5}\sqrt{7},$$
we have shown that $g$ always has real roots, when $s$ and $t$ are distinct.
\medskip

\noindent{(ii) }
Let $m=p=3$ and $\alpha^\DOT=J_2^4,J_1$.
Here, ${\mathcal X}_{J_2,J_2}$ consists of all matrics $X$ of the form
$$
\left[\begin{array}{cccccc}
1&x_{12}&  0   &   0  &   0  &   0  \\
0&  1   &x_{23}&x_{24}&x_{25}&   0  \\
0&   0  &   0  &   0  &   1  &x_{36}\end{array}\right]
$$
and our equations are
$$
\det\zwei{K_3(s)}{X}\ =\ 
\mbox{maximal minors }\zwei{K_2(t)}{X}\ =\ 0
$$
and the same equations with $u$ replacing $t$.
The ideal of these polynomials contains the following univariate polynomial 
$g$, here $e_1=t+u$ and $e_2=tu$.
$$
\begin{array}{c}
x_{36}^3-x_{36}(3s+4e_1)
+x_{36}(4e_1^2+3e_2+10se_1)
-(6e_1e_2+8se_1^2+se_2)\qquad\qquad\qquad\qquad\\
\qquad\qquad\qquad\qquad =\ 
(x_{36}-2e_1)(x_{36}^2-2e_1x_{36}+3e_2)
\ -\ s(x_{36}^2-10e_1x_{36}+8e_1^2+e_2).\rule{0pt}{18pt}
\end{array}
$$
These last two polynomials have roots
$$
e_1\pm \sqrt{e_1^2-3e_2},\  2e_1\quad\mbox{ and }\quad
\frac{5}{3}e_1\pm \frac{\sqrt{e_1^2-3e_2}}{3},
$$
which are interlaced.
For example, if $e_1>0$, then 
$$
e_1- \sqrt{e_1^2-3e_2}<\frac{5}{3}e_1- \frac{\sqrt{e_1^2-3e_2}}{3}<
e_1+ \sqrt{e_1^2-3e_2}<\frac{5}{3}e_1+ \frac{\sqrt{e_1^2-3e_2}}{3}<
2e_1.
$$
When $s,t,u$ are distinct and different from 0, $g$ always has 3 real
roots, by the Intermediate Value Theorem.
We could also note that the discriminant of $g$
$$
\begin{array}{c}
s^2(t-u)^4+t^4(s-u)^2 + u^4(s-t)^2 + s^2t^2(s-t)^2 + s^2u^2(s-u)^2+
\qquad\qquad\qquad\qquad\\
\qquad\qquad (s-t)^2(s-u)^2(t-u)^2 \ +\ 
{\displaystyle \frac{7}{2}}
( s^4(t-u)^2 +t^2(s-u)^4 + u^2(s-t)^4 + t^2u^2(t-u)^2)
\rule{0pt}{18pt}
\end{array}
$$
is a sum of squares and $g(x_{36};1,2,3)$ has 
(floating point) roots
$$
4.736,\  7.756,\  10.508.
$$

\noindent{(iii) }
Let $(m,p)=(3,3)$ and $\alpha^\DOT=(135)^2,J_1^3$.  
Here, ${\mathcal X}_{135,135}$ consists of all matrics $X$ of the form
$$
\left[\begin{array}{ccccccc}
1&x_{12}&   0  &   0  &   0  &   0   \\
0&   0  &   1  &x_{24}&   0  &   0  \\
0&   0  &   0  &   0  &   1  &x_{36}\end{array}\right]
$$
and our equations are
$$
\det\zwei{K_3(s)}{X}\ =\ \det\zwei{K_3(t)}{X}\ =\
\det\zwei{K_3(u)}{X}\ =\ 0.
$$
We write the universal eliminant, $g(x_{36})$, in terms of the elementary
symmetric polynomials in $s,t,u$:
$$
\begin{array}{c}
9x_{36}^6
-48e_1x_{36}^5
+(64e_1^2+108e_2)x_{36}^4
-(288e_1e_2-198e_3)x_{36}^3\qquad\qquad\qquad\\
\qquad\qquad\qquad+(320e_2^2+540e_1e_3)x_{36}^2
-1200e_2e_3x_{36}+1125e_3^2.\rule{0pt}{18pt}
\end{array}
$$ 
Evaluating the parameters $(s,t,u)$ at $(1,2,3)$, we see that
$g(x_{36};1,2,3)$ has roots:
$$
1.491,\ 1.683,\ 3.210,\ 5.630,\ 9.213,\ 10.773.
$$
The discriminant of $g$ is a sum of squares.
The primitive part of the discriminant is
$$
e_3^4
(4e_2^2e_1^2-15e_3e_1^3-15e_2^3+63e_3e_2e_1-81e_3^2)
(256e_2^2e_1^2-768e_3e_1^3-768e_2^3+2592e_3e_2e_1-2187e_3^2)^2.
$$
The second factor is a sum of squares
$$
\frac{7}{2}(s-t)^2(s-u)^2(t-u)^2\ +\ 
\frac{1}{2}s^2((t-u)^4+t^2(s-u)^4+u^2(s-t)^4)
$$
Interestingly, the last (squared) factor is itself a sum of
squares:
$$
\begin{array}{cc}
112(s-t)^2(u^4+s^2t^2)+ 112(t-u)^2(s^4+t^2u^2)
+112(u-s)^2(t^4+s^2u^2)+\qquad\qquad\qquad\\
\qquad 16(s-t)^2(s-u)^2(t-u)^2
+309s^2t^2u^2 
+16(s^4(t^2+u^2)+ t^4(s^2+u^2)+u^4(t^2+u^2)).
\rule{0pt}{18pt}
\end{array}
$$

\noindent{(iv) }
Let $(m,p)=(4,3)$ and $\alpha^\DOT=(135)^4$.  
Here, ${\mathcal X}_{135,135}$ consists of all matrics $X$ of the form
$$
\left[\begin{array}{ccccccc}
1&x_{12}&x_{13}&   0  &   0  &   0  &   0  \\
0&   0  &   1  &x_{24}&x_{25}&   0  &   0  \\
0&   0  &   0  &   0  &   1  &x_{36}&x_{37}\end{array}\right]
$$
and our equations are
$$
\mbox{maximal minors }\zwei{K_3(s)}{X}\ =\ 
\mbox{maximal minors }\zwei{K_5(s)}{X}\ =\ 0
$$
and the same equations with $t$ replacing $s$.
In this case, the universal eliminant has 4 quadratic factors:
$$
\begin{array}{c}
(36x_{12}^2  -x_{12}( 12t+30s)+6st+5s^2), \ \ 
(36x_{12}^2  -x_{12}( 12s +30t) +6st + 5t^2), \\
(3x_{12}^2  -2x_{12}(s+t) +   st), \ \ \mbox{ and }\ \ 
(36x_{12}^2  -30x_{12}(s+t)+ 5t^2+14st+5s^2).\rule{0pt}{15pt}
\end{array}
$$
When $s\neq t$ and neither is zero, we see that each has 2 real roots.
\QED

Observe that in all 4 cases, the discriminant was a sum of squares and the
eliminant has the correct number of real roots for distinct values of the
parameters.
Of particular note is that the system in (ii) was not symmetric in the
parameters and the Schubert data of (iv) was not Pieri Schubert data.
\medskip

Table~\ref{tab:pieri} gives the number of instances of 
Conjecture~\ref{conj:shapiro-grass} we have
checked.

\begin{table}[htb]
 \begin{tabular}{|r||c|c|c||c|c|c|c|c|}
  \hline
  $\alpha^\DOT$:&$(J_2)^5$&$(J_2)^6$&$(J_2)^7$
          &$(J_2)^6$&$(J_3)^5$&$(135)^5$&$(135)^2,(J_1)^6$
             \rule{0pt}{15pt}\\\hline\hline
  $m,p$:& 5,2& 6,2& 7,2&5,3&4,3&5,3&6,3\\\hline
  $d(m,p;\alpha^\DOT)$:&6&15&36&6&16&32&61\\\hline
  \# checked&10,000&2821&504&10160&2002&400&294\\\hline
\end{tabular}\bigskip
 \caption{General Schubert data tested\label{tab:pieri}}
\end{table}

%
%
%
%

\section{Total positivity}\label{sec:tot_pos}
Previous sections have dealt with Schubert conditions given by flags
osculating a real rational normal curve.
Recently, Shapiro and Shapiro have conjectured
that a generalization of this 
choice involving totally positive real matrices would also give only real
solutions. 
We describe that here, prove the first nontrivial instance, and present
some computational evidence in support of this generalization.

A real  upper triangular matrix $g$ with 1's on its diagonal is {\em totally
positive} if every minor of $g$ is positive, except those minors which
vanish on all upper triangular matrices. 
Let $\TP$ be the set of all totally positive,
a multiplicative semigroup.
Define a partial order on real flags $\Fdot$ by $\Fdot<g\Fdot$ if
$g\in\TP$. 

\begin{conj}\label{conj:tot-pos}
For any $m,p>1$, let $\alpha^\DOT$ be
Schubert data for $\mbox{\it Grass}(p,m+p)$.
If $\Fdot^1<\cdots<\Fdot^n$ are real flags, then 
the Schubert varieties $\Omega_{\alpha^1}\Fdot^1,\ldots,
\Omega_{\alpha^n}\Fdot^n$
intersect transversally, with all points of intersection real.
\end{conj}

We will prove Conjecture~\ref{conj:tot-pos} in the first nontrivial case
of $m=p=2$.
First, we relate Conjecture~\ref{conj:tot-pos} to
Conjecture~\ref{conj:shapiro-grass}. 
Let $\Kdot(s)$ be the square matrix of size $(m+p)$ whose $i,j$th entry is
$\binom{j-i}{i-1}s^{j-i}$ ({\it cf.}~(\ref{eq:K-matrix})).
If $s>0$, then $\Kdot(s)$ is totally positive
and for any $s,t$ we have $\Kdot(s)\cdot\Kdot(t)=\Kdot(s+t)$.
To see this, first recall that $\TP$ is generated as a semigroup by
$\exp(E_{i,i+1})$, where $E_{i,i+1}$ is the elementary matrix whose only
non-zero entry is in position $i,i+1$~\cite{Loewner}.
These assertions follow from the observation that
$$
\Kdot(s)\ =\ \exp(s N),
$$
where $N$ is the nilpotent matrix whose only non-zero entries are
$(1,2,\ldots,m+p-1)$ lying just above its main diagonal.

Theorem~\ref{thm:implies} holds in
this new setting.
For this, we alter the notion of Pieri Schubert data $\alpha^\DOT$
to Schubert data $\alpha^1,\ldots,\alpha^n$ where all except
possibly $\alpha^1$ and $\alpha^n$ are Pieri conditions.

\begin{thm}\label{thm:tot-pos}
Let $a,b>1$ and suppose that Conjecture~\ref{conj:tot-pos} holds for 
$(m,p)=(a,b)$ and Schubert data $\alpha^\DOT=(J_1)^{mp}$.
Then Conjecture~\ref{conj:tot-pos} holds for any Pieri Schubert data for
$\mbox{\it Grass}(p,m+p)$ where $(m,p)\leq (a,b)$ or $(b,a)$
coordinatewise. 
\end{thm}

\noindent{\bf Proof. }
The arguments used to prove Theorem~\ref{thm:implies} work 
here with minor adjustments.

We first remark that total positivity, and hence our order $<$ on real
flags, is defined with respect to a choice of ordered basis for 
${\mathbb R}^{m+p}$.
Suppose that $e_1,\ldots,e_{m+p}$ is the basis we used to define this order.
Then $\Fdot<G_\DOT$ is and only if $G_\DOT<'\Fdot$, where $<'$ is defined
with respect to the basis $e_1,-e_2,e_3,-e_4,\ldots$.
Similarly, if the basis $e_1,\ldots,e_{m+p}$ is orthonormal, 
then $\Fdot<G_\DOT$ if and only if 
$F_\DOT^\perp<''G_\DOT^\perp$, where $<''$ is defined with respect to
the basis in reverse order $e_{m+p},\ldots,e_2,e_1$.
Thus 
$$
F_\DOT^1<F_\DOT^2<\cdots<F_\DOT^n\ 
\Longleftrightarrow\ 
F_\DOT^n<'\cdots<'F_\DOT^2<'F_\DOT^1
$$
so that Conjecture~\ref{conj:tot-pos} holds for Schubert data $\alpha^\DOT$
if and only if it holds for the data in reverse order.
(This is the only rearrangment we used in the proof of
Lemma~\ref{lem:ordering}.) 
Similarly, the analogue of Lemma~\ref{lem:ordering}(i) holds.
For  the analogue of  Lemma~\ref{lem:ordering}(ii),
permute the last two Schubert conditions, so that $\beta^\DOT$ is
still Pieri Schubert data, in our new, restricted definition.

Finally, in the proof of Lemma~\ref{lem:ordering}(iii), replace
$s_3,\ldots,s_n$ in defining the set ${\mathcal O}$ by fixing $F_\DOT^1$ to 
be the standard flag represented by the matrix $I_{m+p}$ and 
let ${\mathcal O}$ be the set of all
$$
F_\DOT^1<\cdots<F_\DOT^n
$$
where the appropriate transversality conditions hold.
Since $\TP$ is open, it follows that there exists $\epsilon>0$ and 
totally positive matrix $M$ (which stabilizes $F_\DOT^1$) such that 
if $0<s<\epsilon$, then 
$F_\DOT^1< M\cdot\Kdot(s)\cdot F_\DOT^1< F_\DOT^2$.
Then the same arguments used to prove Theorem~\ref{thm:implies}
suffice.
In particular, the analog of Proposition~\ref{no-mult} also holds in this
setting.
\QED

Totally positive matrices have a useful description. 
Let ${\mathcal U}$ be the group of real unipotent (upper triangular)
matrices. 
Then $\TP$ is a connected component of the complement of
a hypersurface $H$ ${\mathcal U}$ defined by the vanishing of all minors
consisting of the first $i$ rows and last $i$
columns~\cite{Shapiro_TP}. 
This has a geometric description.

Associating a matrix to a flag as in Section~\ref{sec:loc-coords}, we may
identify ${\mathcal U}$ with a Zariski open subset of 
the real flag manifold.
Then the hypersurface $H$ is the union of all positive
codimension Schubert varieties defined by the flag determined by the
identity matrix.

Given a matrix $M\in{\mathcal U}$, the translate
$\TP.M$ is a component of the complement of all Schubert varieties of
positive codimension defined by the flag given by $M$. 
Similarly, given a totally positive matrix $M$, the set of upper triangular
matrices  $N$ for which there exists a totally positive $g$ with $gN=M$ is
the component of this complement containing the
identity matrix.

Let $F_\DOT^1<\cdots<F_\DOT^n$ be real flags.
Using a real automorphism of the flag manifold, we may assume that
$F_\DOT^1=\Kdot(0)=I_{m+p}$.
Then $F_\DOT^2,\ldots,F_\DOT^n\in\TP$, since they are all
translates of the identity by totally positive matrices.
Also,  $F_\DOT^1,\ldots,F_\DOT^{n-1}$ are in the same component of the
complement  
of all positive dimensional Schubert cells defined by $F_\DOT^n$.
If we now consider a real coordinate transformation fixing $F_\DOT^1$, but
with $F_\DOT^n$ becoming $\Kdot(\infty)$, then this complement becomes 
$\TP$, in these new coordinates.

Thus we may work in the local coordinates
${\mathcal X}:={\mathcal X}_{\alpha^1,\alpha^n}$.
We do this in our proof of the following theorem
and in subsequent calculations.

\begin{thm}\label{thm:tp22}
Conjecture~\ref{conj:tot-pos} holds for $m=p=2$ and Schubert data
$(J_1)^4$.
\end{thm}

\noindent{\bf Proof. }
Let $F,G\in \TP$ be totally positive matrices
and set $H = G\cdot F$.
When $m=p=2$, ${\mathcal X}={\mathcal X}_{J_1,J_1}$ is the set of matrices 
$$
\left[
\begin{array}{cccc}1&a&0&0\\0&0&1&b\end{array}\right].
$$

For a matrix $L$, let 
$L_{ij}$ denote the $2\times 2$-minor of $L$ given by the first two rows
and columns $i$ and $j$.
Then the equations for a 2-plane in ${\mathcal X}$ to meet the flags given
by $F$ and $H$ are
\begin{eqnarray*}
f&:=& F_{24} - b F_{23} - a F_{14} + ab F_{13}\\
h&:=& H_{24} - b H_{23} - a H_{14} + ab H_{13}
\end{eqnarray*}
The lexicographic Gr\"obner basis for this (with $a<b$) is
\begin{eqnarray*}
H_{13}f-F_{13}h&=&   J_{14} - b J_{24} -a J_{34}\\
(H_{14}-b H_{13})f-(F_{14}-b F_{13})h
&=&J_{13}-b(J_{23}+J_{14}) + b^2 J_{24},
\end{eqnarray*}
where $J_{ij}$ is the $ij$th minor of the matrix:
$$
\left[
\begin{array}{cccc}
F_{24}& F_{23}& F_{14}&  F_{13}\\
H_{24}& H_{23}& H_{14}&  H_{13}\end{array}\right].
$$
We may write the the discriminant of the quadratic
equation for $b$ as follows
$$
(J_{23}+J_{14})^2-4J_{13}J_{24}\ =
\ (L_{23}+L_{14})^2-4L_{13}L_{24},
$$
where $L$ is the matrix
$$
\left[
\begin{array}{cccc}
F_{13}& F_{14}& H_{13}& H_{14}\\
F_{23}& F_{24}& H_{23}& H_{24}\end{array}\right].
$$
Thus we will have two real roots for our original system if and only if 
$$
\Lambda(B)\ :=\ L_{13}-B(L_{23}+L_{14}) + B^2 L_{24}\ =\ 0
$$
has 2 real solutions.
Painstaking calculations reveal that 
$\Lambda(1)=-G_{12}G_{34}<0$.
Since $L_{24}=H_{13}H_{24}-H_{23}H_{14} = H_{12}H_{34}$
by the Pl\"ucker relations, we see that 
$L_{24}>0$ and so $\Lambda(B)=0$ will have 2 real solutions.
\QED

Table~\ref{table:tpinstances} shows the number of instances of
Conjecture~\ref{conj:tot-pos} that we have verified.

\begin{table}[htb]
 \begin{tabular}{|c||c|c|c|c|c|c|}
  \hline $\alpha^\DOT$&$(J_1)^6$&$(J_2)^5$&$(135)^4$&$(J_1)^8$
     &$(J_2)^6$&$(135)(136)(J_1)^5$\rule{0pt}{12pt}\\\hline\hline
  $m,p$&3,2&5,2&4,3&4,2&6,2&4,3\\\hline
  $d$&5&6&8&14&15&25\\\hline
  \# checked&12000&4000&4000&1500&300&150\\\hline
 \end{tabular}\bigskip
 \caption{Instances checked\label{table:tpinstances}}
\end{table}
 
\section{Further remarks}\label{sec:counter}
We present a counterexample to the original conjecture of Shapiro and
Shapiro and close
with a discussion of further questions.

\subsection{A counterexample to the original conjecture}
The original conjecture of Shapiro and Shapiro concerned the $M$-property
for flag manifolds~\cite{Shapiro_M}.
An algebraic set $X$ defined over ${\mathbb R}$ has the $M$-property if
the sum of the ${\mathbb Z}/2{\mathbb Z}$-Betti numbers of $X({\mathbb R})$
and of $X({\mathbb C})$ are equal.
Shapiro and Shapiro conjectured that an intersection of Schubert cells in a
flag manifold has the $M$-property, if the cells are defined by flags
osculating the rational normal curve at real points.
When such an intersection is zero-dimensional all of
its points are real.
It is this consequence we have been studying.

While there is much evidence in support of this conjecture for zero
dimensional intersections in a Grassmannian
(Conjectures~\ref{conj:shapiroII},~\ref{conj:shapiro-grass},
and~\ref{conj:tot-pos}),
it does not hold for more general flag manifolds.
In fact, we give a counter example in the simplest enumerative problem in a
flag manifold that does not reduce to an enumerative problem in
a Grassmannian. 

\begin{cex}\label{cex}
Consider the manifold ${\mathbb F}(2,3;5)$ consisting of partial 
flags $X\subset Y$
in ${\mathbb C}^5$ with $\dim X=2$ and $\dim Y=3$.
This manifold has dimension 8;
the projection to $\mbox{\it Grass}(2,5)$ has fibre over a 2-plane $X$ equal
to ${\mathbb P}({\mathbb C}^5/X)\simeq{\mathbb P}^2$.
Given general 2-planes $a,b$, and $c$ and general 3-planes $A,B$, and $C$,
there are 4 flags $X\subset Y$ which satisfy
\begin{equation}\label{eq:flag-cond}
 \begin{array}{c}
  \mbox{$X$ meets $a,B$, and $C$ nontrivially}\\
  \mbox{$\dim Y\cap A\geq 2$ and $Y$ meets $b$ and $c$ nontrivially}
 \end{array}
\end{equation}
That this number is 4 may be verified using the Schubert calculus for a flag
manifold~\cite{Fulton_tableaux} or the equations we give below.

Let $\Kdot(s)$ be the flag of subspaces osculating the standard rational
normal curve.
Set
$$
\begin{array}{ccc}
a\ :=\ K_2(4)&\quad& A\ :=K_3(0)\\ 
b\ :=\ K_2(1)&     & B\ :=K_3(3)\\ 
c\ :=\ K_2(-5)&    & C\ :=K_3(-1)
\end{array}
$$
We claim that of the 4 flags $X\subset Y$ satisfying~(\ref{eq:flag-cond})
for this choice of $a,b,c,A,B,C$, 2 are real and 2 are complex.

We outline the computation.
Choose local coordinates for ${\mathbb F}(2,3;5)$ as follows.
Let $Y$ be the row space of the $3\times 5$-matrix
$$
\left[\begin{array}{ccccc}
0&0&  1   &x_{14}&x_{15}\\
1&0&x_{23}&x_{24}&x_{25}\\
0&1&x_{33}&x_{34}&x_{35}\end{array}\right]
$$
and $X$ be the row space of its last 2 rows.
We seek the solutions to the following overdetermined system
of polynomials: 
$$
\begin{array}{c}
\det\zwei{K_2(1)}{Y}\ =\ \det\zwei{K_2(-5)}{Y}\ =\ 
\det\zwei{K_3(3)}{X}\ =\ \det\zwei{K_3(-1)}{X}\ = \\ \rule{0pt}{25pt}
\mbox{maximal minors }\zwei{K_2(4)}{X}\ =\ 
\mbox{maximal minors }\zwei{K_3(0)}{Y}\ =\ 0
\end{array}
$$
These polynomials generate a zero-dimensional ideal 
containing the following univariate polynomial, which is part of a
lexicographic Gr\"obner basis satisfying the Shape Lemma:
$$
27063 - 117556 x_{14} - 5952 x_{14}^2
-10416 x_{14}^3 + 32400 x_{14}^4.
$$
This has roots
$$
-.736 \pm 1.30 \sqrt{-1}, \quad .227, \quad 1.62.\
$$
Thus 2 of the flags are complex.
\end{cex}

\subsection{Further questions}

While Counterexample~\ref{cex} shows that we cannot guarantee all points of
intersection real when the Schubert varieties are given by flags osculating
a real rational normal curve, a number of questions remain (besides the
resolution of the conjectures of the previous sections).
There remains the original question of Fulton:
\medskip

\noindent{\bf Question 1: }
Given Schubert data for a flag manifold, do there exist real flags in general
position whose corresponding Schubert varieties have {\it only} real points
of intersection? 
\medskip

In every case we know, this does happen.
For instance, if we change the 3-plane $B$ to $K_3(2)$ in
Counterexample~\ref{cex}, then all 4 solution flags are real. 
There is also the following result, showing this holds in infinitely many
cases.
Let ${\mathbb F}(2,n-2;n)$ be the manifold of flags $X\subset Y$
in ${\mathbb C}^n$ where $\dim X=2$ and $\dim Y=n-2$.
A Grassmannian Schubert condition is a Schubert condition on the flag 
$X\subset Y$  which only imposes conditions on one of the subspaces.
We likewise define Grassmannian Schubert data.
For example, Counterexample~\ref{cex} involves 
Grassmannian Schubert data.

\begin{prop}[Theorem 13 of~\cite{Sottile97b}]\label{prop:MEGA}
Given any Grassmannian Schubert data for ${\mathbb F}(2,n-2;n)$,
there exist real flags whose corresponding Schubert varieties meet
transversally with all points of intersection real.
\end{prop}

The beauty of the conjectures of Shapiro and Shapiro is that they give a
simple algorithm for selecting the flags defining the Schubert varieties.
\medskip

\noindent{\bf Question 2: }
Can the choice of flags in Question 1 (or Proposition~\ref{prop:MEGA})
be made effective?
In particular, is there an algorithm for selecting these flags?
\medskip

While computing the examples described here,
we have made a number of observations which deserve further scrutiny.
These concern eliminant polynomials in the ideals defining the intersections
of Schubert varieties in the local coordinates we have been using.

Suppose we have Schubert data $\alpha^\DOT$, and have chosen local
coordinates either for the Grassmannian or are working in 
${\mathcal X}_{\alpha^n,\alpha^{n-1}}$.
Conjecture~\ref{conj:shapiro-grass} or~\ref{conj:tot-pos} may be formulated
in terms of a parameterized system of polynomials with parameters either
$s_1,\ldots,s_n$ in the case of Conjecture~\ref{conj:shapiro-grass} or 
$(n-1)$-tuples of totally positive matrices (or in terms of some
parameterization of $\TP$~\cite{BFZ}).
For each of the coordinates, the ideal of this system contains a
universal eliminant, which is the minimal univariate polynomial in that
coordinate with coefficients rational functions in the parameters.

We ask the following questions about the eliminant.
\medskip

\noindent{\bf Question 3: }
Does the universal eliminant have degree equal to the generic number of
solutions?
That is, do generic solutions satisfy the shape lemma?
\medskip

\noindent{\bf Question 4: }
Let $\Delta$ be the discriminant of the polynomial system, 
a polynomial in the parameters which vanishes when there are solutions with
multiplicities. 
\begin{enumerate}
 \item[a)]
   Is the locus $\Delta\neq 0$ connected?

 \item[b)]
   In the case of Conjecture~\ref{conj:shapiro-grass}, where $\Delta$ is
   a polynomial in the parameters $s_1,\ldots,s_n$, is $\Delta$ always a
   sum of squares of polynomials?
 
 \item[c)]
   If so, are these polynomials monomials in the $s_i$ and their
   differences $(s_i-s_j)$? 
   This would imply that the polynomial systems are always multiplicity-free
   for distinct real values of the parameters, and hence the stronger
   version of Theorem~\ref{thm:implies} mentioned in
   Remark~\ref{rem:stronger}.
\end{enumerate}
\medskip

The discriminants we have computed for instances of the conjectures
for the Grassmannian (including the discriminant for system of
Theorem~\ref{thm:tp22}) are always non-negative when the parameters are
distinct.
For the case of Counterexample~\ref{cex}, we computed a
discriminant for a simpler, but equivalent system, in the spirit of
sections~\ref{sec:simpler} and~\ref{sec:loc-coords}.
This polynomial in parameters $s_1,s_2,t_1,t_2$ is symmetric in
the $s$'s and in the $t$'s separately (and in the transformation
$s_i\leftrightarrow t_i$) and has degree 24.
It has three factors, the first of degree 20 with 857 terms,
and the square
$$
(2s_1s_2+2t_1t_2-(s_1+s_2)(t_1+t_2))^2.
$$
While this factor will not prevent the discriminant from being a sum of
squares, this factor shows that there is
a choice of distinct parameters for which the discriminant vanishes.
Indeed, if we set $s_1=3,s_2=6,t_1=9$, and $t_2=5$, then this factor
vanishes, and the resulting system has a root of multiplicity 2.
This also explains why different values of the parameters in
Counterexample~\ref{cex} give different numbers of real and complex
solutions.
\medskip

\noindent{\bf Question 5: }
When the universal eliminant factors over ${\mathbb Z}$, it 
reflects either some underlying geometry or some interesting arithmetic.
More generally, one might ask about the Galois group of these enumerative
problems~\cite{Harris_galois}, 
or the Galois group of the universal eliminant.
For instance, is it the full symmetric group?
That is not always the case, as the example
of Theorem~\ref{thm:grass-proof}(iv) shows.
\medskip

\noindent{\bf Question 6: }
In many cases with the substitution of $s_i=i$,
the eliminant factors over the integers. 
This happens in Conjecture~\ref{conj:shapiroI}, Theorem~\ref{thm:23},
Theorem~\ref{thm:grass-proof}(i) and (iv), and in other cases.
Table~\ref{table:factor} lists the degrees of the factors in the case of 
Conjecture~\ref{conj:shapiroI}.
\begin{table}[htb]
 \begin{tabular}{|c||c|c|c|c|c|c|c|}
  \hline$m,p$&3,2&4,2&5,2&6,2&7,2&3,3&3,4\\\hline\hline
  $d_{m,p}$&5&14&42&132&429&42&462\\\hline
  Factors&2,3&6,8&10,32&20,112&&6,36&16,30,416\\\hline
 \end{tabular}\bigskip
 \caption{Factorization of the eliminant\label{table:factor}}
\end{table}
 Why does this choice of $s_i=i$ induce a factorization?
 Is there any special geometry or interesting arithmetic here?
 If 2 parameters are allowed to come together, then 
 the resulting ideal factors in a way respecting the
 product of Schubert classes, by the Corollary to Theorem 1
 in~\cite{EH87}.
 From the Schubert calculus, we would expect factors of 9 and 5 for
 $(m,p)=(2,4)$, 14 and 28 for $(m,p)=(2,5)$, and 21 and 21 for 
 $(m,p)=(3,3)$, but these do not appear in Table~\ref{table:factor}.

\section*{Acknowledgements}
I thank Boris Shapiro for sharing these conjectures 
and Joachim Rosenthal with whom I found the first strong evidence for their
validity~\cite{RS98}.
Bernd Sturmfels helped me learn computer algebra and
suggested using discriminants.
This work also benefitted from discussions with 
David Eisenbud,  Ioannis Emiris, Jean-Charles Faug\`ere, Bill
Fulton, Birk Huber,  Bernard Mourrain, Bruce Reznick,
Fabrice Rouillier, Jan Verschelde,
and Paul Zimmermann. 
Lastly, this would not be possible were it not
for the authors of the following software packages:
Linux, Maple, Singular, Macaulay2, MuPAD, FGB, RealSolving, PHC,
and also my brothers Larry Sottile and Joe Sottile who helped me build the
computer which performed most of these experiments.

\end{document}